\documentclass[11pt]{article}
\usepackage{amssymb}
\usepackage{graphicx}

\newcommand{\nuc}{\newcommand}

\nuc{\plim}{{_p \atop ^{^{\ra}}}}
\nuc{\sgg}{\sg^2}

\nuc{\bs}{\begin{slide*}{}}
\nuc{\es}{\end{slide*}}
\nuc{\bc}{\begin{center}}
\nuc{\ec}{\end{center}}
\nuc{\bul}{\bullet}
\nuc{\ze}{\zeta}

\nuc{\bv}{\begin{verbatim}}
\nuc{\ev}{\end{verbatim}}

\nuc{\bt}{\begin{tabular}}
\nuc{\et}{\end{tabular}}
\nuc{\lera}{\leftrightarrow}

\nuc{\BT}{\begin{tabbing}}
\nuc{\ET}{\end{tabbing}}

\newcommand{\tht}{\theta}
\nuc{\Th}{\Theta}
\newcommand{\gam}{\gamma}
\nuc{\Gam}{\Gamma}
\newcommand{\om}{\omega}
\nuc{\Om}{\Omega}
\newcommand{\al}{\alpha}

\newcommand{\be}{\beta}

\nuc{\pde}{\partial}
\newcommand{\lam}{\lambda}

\nuc{\ph}{\phi}
\nuc{\Ph}{\Phi}
\nuc{\zi}{\xi}
\nuc{\Zi}{\Xi}
\nuc{\ep}{\epsilon}
\nuc{\Ep}{\Epsilon}
\nuc{\del}{\delta}
\nuc{\Del}{\Delta}
\nuc{\sg}{\sigma}
\nuc{\Sg}{\Sigma}
\nuc{\pa}{\parallel}
\nuc{\DB}{\begin{description}}
\nuc{\DE}{\end{description}}
\nuc{\cc}{\cite}
\nuc{\lab}{\label}
\nuc{\half}{\frac{1}{2}}
\nuc{\ha}{\frac{1}{2}}
\nuc{\hl}{\hline}
\nuc{\sco}{\setcounter{equation}{0}}

\nuc{\EQ}{\begin{displaymath}}
\nuc{\EN}{\end{displaymath}}
\newcommand{\EQN}{\begin{equation}}
\newcommand{\ENN}{\end{equation}}
\nuc{\EBN}{\begin{eqnarray}}
\nuc{\EEN}{\end{eqnarray}}
\newcommand{\EB}{\begin{eqnarray*}}
\newcommand{\EE}{\end{eqnarray*}}
\nuc{\ol}{\overline}
\nuc{\ul}{\underline}

\nuc{\nn}{\nonumber}
\nuc{\Ra}{\Rightarrow}
\nuc{\ra}{\rightarrow}
\nuc{\LRA}{\Leftrightarrow}
\nuc{\Lra}{\Longrightarrow}
\nuc{\lra}{\longrightarrow}
\nuc{\lba}{\left(\begin{array}}
\nuc{\lcba}{\left\{\begin{array}}
\nuc{\ecar}{\end{array}\right)}
\nuc{\ear}{\end{array}\right)}
\nuc{\bnar}{\begin{array}}
\nuc{\enar}{\end{array}}
\nuc{\til}{\tilde}
\nuc{\stk}{\stackrel}

\nuc{\bay}{\begin{array}}
\nuc{\eay}{\end{array}}

\nuc{\noi}{\noindent}
\nuc{\qm}{q^{-1}}
\nuc{\zm}{z^{-1}}

\newcounter{expl}
\nuc{\atce}{\addtocounter{expl}{1}}
\nuc{\abce}{\arabic{expl}}
\nuc{\abcc}{\arabic{chapter}}

\newcounter{taybl}
\nuc{\atct}{\addtocounter{taybl}{1}}
\nuc{\abct}{\arabic{taybl}}

\nuc{\ben}{\begin{enumerate}}
\nuc{\een}{\end{enumerate}}

\newcounter{vsec}
\nuc{\vstep}{\stepcounter{vsec}}
\nuc{\vreset}{\setcounter{section}{\value{vsec}}}

\nuc{\scone}{\setcounter{section}{1}}
\nuc{\sctwo}{\setcounter{section}{2}}
\nuc{\scthree}{\setcounter{section}{3}}
\nuc{\scfour}{\setcounter{section}{4}}
\nuc{\scfive}{\setcounter{section}{5}}
\nuc{\scsix}{\setcounter{section}{6}}
\nuc{\scseven}{\setcounter{section}{7}}
\nuc{\sceight}{\setcounter{section}{8}}
\nuc{\scnine}{\setcounter{section}{9}}
\nuc{\scten}{\setcounter{section}{10}}
\nuc{\sceleven}{\setcounter{section}{11}}
\nuc{\sctwelve}{\setcounter{section}{12}}
\nuc{\scthirteen}{\setcounter{section}{13}}
\nuc{\scfourteen}{\setcounter{section}{14}}
\nuc{\scfifteen}{\setcounter{section}{15}}
\nuc{\scsixteen}{\setcounter{section}{16}}
\nuc{\scseventeen}{\setcounter{section}{17}}
\nuc{\sceighteen}{\setcounter{section}{18}}
\nuc{\scnineteen}{\setcounter{section}{19}}
\nuc{\sctwenty}{\setcounter{section}{20}}

\nuc{\bit}{\hspace{.1em}}

\nuc{\mbo}{$\mbox{ }$}
\nuc{\mbh}{\mbo\mbo\mbo}
\nuc{\mucl}{\multicolumn}
\nuc{\renu}{\renewcommand}
\nuc{\swdth}{\settowidth}
\nuc{\awdth}{\addtolength}
\newlength{\wdth}

\newlength{\wdths}
\nuc{\hspw}{\hspace{\wdth}}
\nuc{\hspws}{\hspace{\wdths}}
\nuc{\swid}[1]{\swdth{\wdth}{#1}}
\nuc{\ens}{\ensuremath}
\nuc{\nors}{\normalsize}
\nuc{\benq}{\ben}
\nuc{\enq}{\een}
\nuc{\itmq}{\item[] }
\nuc{\itemq}{\item[] }

\nuc{\setc}{\setcounter}
\nuc{\stepc}{\stepcounter}
\nuc{\nuco}{\newcounter}
\nuc{\rstepc}{\refstepcounter}

\nuc{\remk}[1]{\ens{\ul{Remarks} \mbo\kstp{#1.\kval}. }}
\nuc{\vchap}{\arabic{chapter}}

\nuc{\subst}{\subsection*}
\nuc{\ssubst}{\subsubsection*}
\nuc{\secst}{\section*}
\nuc{\sx}[1]{\secst{\vstp{\vval.} #1}
\setc{xtem}{0}\setc{etem}{0}\setc{ttem}{0}\setc{dtem}{0}}
\nuc{\atocz}[1]{\addcontentsline{toc}{section}{\vval. \hspace{.3em}#1}}
\nuc{\sectoc}[1]{\sx{#1}\atocz{#1}} 
\nuc{\atoc}[1]{\addcontentsline{toc}{subsection}{#1}}
\nuc{\atocs}[1]{\addcontentsline{toc}{subsection}{\small\it #1 \nors}}
\nuc{\atocx}[1]{\addcontentsline{toc}{subsection}{\small\it\vchap.\val. #1 \nors}}
\nuc{\subtoc}[1]{\subsection*{#1}\atocs{#1}}

\nuc{\aob}[2]{\ens\big({#1 \atop #2}\big)}
\nuc{\aobc}[2]{\ens\big\{{#1 \atop #2}\big\}}
\nuc{\aonb}[2]{\ens\frac{#1}{#2}}
\nuc{\aoboc}[3]{\lba{c}#1\\#2\\#3\ear}
\nuc{\aobcod}[4]{\ens\big({#1 \atop #2}{\mbox{ } #3 \atop \mbox{ } #4}\big)}
\nuc{\aobcodx}[4]{\ens\left[{#1 \atop #2}{\mbox{ } #3 \atop \mbox{ } #4}\right]}
\nuc{\aobcodL}[4]{\ens\left[{#1 \atop #2}{ #3 \atop  #4}\right]}
\nuc{\abx}[2]{ {{#1 \atop #2}}}
\nuc{\abL}[2]{\left[ {#1 \atop #2}\right]}
\nuc{\abc}[3]{ \left[ { { #1 \atop #2} \atop _{#3}}\right]}
\nuc{\abcL}[3]{ \left[ { { #1 \atop #2} \atop _{#3}}\right]}
\nuc{\abcx}[3]{  { { #1 \atop #2} \atop _{#3}}}
\nuc{\vecL}[3]{ \left[ { { #1 \atop #2} \atop _{#3}}\right]}
\nuc{\vecx}[3]{  { { #1 \atop #2} \atop _{#3}}}
\nuc{\abcdx}[4]{{ { #1 \atop #2} \atop {#3 \atop #4}}}
\nuc{\ivecx}[4]{{ { #1 \atop #2} \atop {#3 \atop #4}}}
\nuc{\abcdL}[4]{\left[{ { #1 \atop #2} \atop {#3 \atop #4}}\right]}
\nuc{\ivecL}[4]{\left[{ { #1 \atop #2} \atop {#3 \atop #4}}\right]}
\nuc{\tubytu}[4]{ \left[ \abx{#1}{#2} \abx{#3}{#4} \right]}
\nuc{\abcS}[6]{ \left[ \abcx{#1}{#2}{#3} \abcx{#4}{#5}{#6} \right]}
\nuc{\thbytu}[6]{ \left[ \abcx{#1}{#2}{#3} \abcx{#4}{#5}{#6} \right]}
\nuc{\abcs}[6]{\left[ { #1 \atop #2} { #3 \atop #4} { #5 \atop #6} \right]}
\nuc{\tubyth}[6]{\left[ { #1 \atop #2} { #3 \atop #4} { #5 \atop #6} \right]}
\nuc{\abcM}[9]{ \left[ \abcx{#1}{#2}{#3} \abcx{#4}{#5}{#6} 
\abcx{#7}{#8}{#9}\right]}
\nuc{\thbyth}[9]{ \left[ \abcx{#1}{#2}{#3} \abcx{#4}{#5}{#6} 
\abcx{#7}{#8}{#9}\right]}
\nuc{\abcdex}[5]{{ {{{ #1 \atop #2} \atop {#3 \atop #4}} \atop _{#5}}}  }
\nuc{\abcdeL}[5]
{\left[{ {{{ #1 \atop #2} \atop {#3 \atop #4}} \atop _{#5}}}\right]}
\nuc{\vvecx}[5]{{ {{{ #1 \atop #2} \atop {#3 \atop #4}} \atop _{#5}}}  }
\nuc{\vvecL}[5]
{\left[{ {{{ #1 \atop #2} \atop {#3 \atop #4}} \atop _{#5}}}\right]}
\nuc{\dadb}[2]{\ens{\frac{\pde #1}{\pde #2}}}
\nuc{\dadbd}[2]{\ens{\frac{d #1}{d #2}}}
\nuc{\dsqdadb}[2]{\ens{\frac{\pde^2}{\pde #1 \pde #2}}}
\nuc{\dsqadbdc}[3]{\ens{\frac{\pde^2 #1}{\pde #2 \pde #3}}}
\nuc{\dsqadbdcd}[3]{\ens{\frac{d^2 #1}{d #2 d #3}}}
\nuc{\dsqadbsq}[2]{\ens{\frac{\pde^2 #1}{\pde #2^2}}}
\nuc{\dsqadbsqd}[2]{\ens{\frac{d^2 #1}{d #2^2}}}
\nuc{\adotb}[2]{\ens{#1,\cdots,#2}}
\nuc{\dash}{$"$}
\nuc{\dai}{^{'}}
\nuc{\dae}{$^{'}$}
\nuc{\daes}{$^{'}$ }
\nuc{\ddai}{^{''}}
\nuc{\ddae}{$^{''}$}
\nuc{\emu}[1]{\ens{e^{-#1}}}
\nuc{\epu}[1]{\ens{e^{#1}}}
\nuc{\asubb}[2]{\ens{#1_{#2}}}
\nuc{\calx}[1]{\ens{{\cal #1}}}
\nuc{\subc}[1]{\ens{_{(#1)}}}
\nuc{\aorb}[5]{ #1 = \lcba {r@{\quad if\quad}l} #2 & #3 \\ #4 & #5 \end{array}\right. }
\nuc{\aorborc}[7]{ #1 = \lcba {r@{\quad:\quad}l} #2 & #3 \\ #4 & #5\\ #6 & #7 \end{array}\right. }
\nuc{\aorbx}[6]{ #1 #6 \lcba {r@{\quad if\quad}l} #2 & #3 \\ #4 & #5 \end{array}\right. }
\nuc{\aorbs}[4]{\lcba {r@{\quad if\quad}l} #1 & #2 \\ #3 & #4 \end{array}\right. }

\nuc{\intab}[2]{\ens{\int_{#1}^{#2}}}

\nuc{\mbspace}[1]{\ens{$\mbox{ }$\hspace{#1 em}}}

\nuc{\subab}[2]{\ens{#1_#2}}  
\nuc{\wov}[1]{\ens{\frac{1}{#1}}} 
\nuc{\sumab}[2]{\ens{\Sg_#1^#2}}  
\nuc{\wtv}[1]{\ens{1,\cdots,#1}} 
\nuc{\vtv}[2]{\ens{#1,\cdots,#2}} 
\nuc{\upcx}[1]{\ens{^{(#1)}}}  
\nuc{\subcx}[1]{\ens{_{(#1)}}}  
\nuc{\und}[1]{_{_#1}}
\nuc{\snk}[1]{#1_{n,k}}   

\nuc{\wfig}{\begin{wrapfigure}{r}{0.5\textwidth}}
\nuc{\winc}[1]{\includegraphics[width=0.48\textwidth]{#1}}
\nuc{\wfigb}[1]{\begin{wrapfigure}[#1]{r}{0.5\textwidth}}
\nuc{\wfige}{\end{wrapfigure}}

\nuc{\alr}{\al_r}

\nuc{\bmp}{\begin{minipage}}
\nuc{\bfig}{\begin{figure}}
\nuc{\bpf}{\begin{proof}}

\nuc{\cip}{{_{p} \atop ^{\ra}}}
\nuc{\call}{{\cal L}}
\nuc{\calh}{{\cal H}}
\nuc{\cale}{{\cal E}}
\nuc{\cali}{{\cal I}}
\nuc{\calp}{{\cal P}}
\nuc{\calf}{{\cal F}}
\nuc{\calg}{{\cal G}}
\nuc{\cald}{{\cal D}}

\nuc{\colw}{\columnwidth}

\nuc{\domtpi}{\frac{d\om}{2\pi}}

\nuc{\ept}{\ep_t}
\nuc{\epi}{\ep_i}
\nuc{\emp}{\end{minipage}}

\nuc{\efig}{\end{figure}}
\nuc{\epf}{\end{proof}}

\nuc{\fh}{\hat{f}}

\nuc{\intoT}{\int_0^T}
\nuc{\intot}{\int_0^t}

\nuc{\intoi}{\int_0^\infty}
\nuc{\intmii}{\int_{-\infty}^\infty}
\nuc{\intow}{\int_0^1}
\nuc{\intmww}{\int_{-1}^1}
\nuc{\intmwo}{\int_{-1}^0}
\nuc{\ion}{\frac{i}{n}}
\nuc{\intiin}{\int_{\frac{i-1}{n}}^{\frac{i}{n}}}
\nuc{\intmpp}{\int_{-\pi}^{\pi}}
\nuc{\intopiot}{\int_0^{\pi/2}}

\nuc{\limni}{\lim_{n\ra\infty}}
\nuc{\limepo}{\lim_{\ep\ra 0}}
\nuc{\limTi}{\lim_{T\rai}}

\nuc{\nut}{\nu_t}

\nuc{\omk}{\om_k}
\nuc{\odel}{o(\del)}

\nuc{\phir}{\phi_r}

\nuc{\quart}{\frac{1}{4}}
\nuc{\qh}{\hat{q}}
\nuc{\Qh}{\hat{Q}}

\nuc{\rai}{\ra\infty}
\nuc{\rao}{\ra 0}

\nuc{\sgwn}{\Sg_1^n}
\nuc{\sgwp}{\Sg_1^p}
\nuc{\sgwm}{\Sg_1^m}
\nuc{\sgmmm}{\Sg_{-m}^m}
\nuc{\sgwi}{\Sg_1^\infty}
\nuc{\sgmii}{\Sg_{-\infty}^\infty}
\nuc{\sgoi}{\Sg_0^\infty}
\nuc{\sgwN}{\Sg_1^N}
\nuc{\sgwr}{\Sg_1^r}
\nuc{\sgwM}{\Sg_1^M}

\nuc{\sgwT}{\Sg_1^T}

\nuc{\tpikon}{\frac{2\pi}{n}}

\nuc{\upm}{^{-1}}
\nuc{\upgi}{^{-}}
\nuc{\upa}{\uparrow}

\nuc{\woT}{\frac{1}{T}}
\nuc{\won}{\frac{1}{n}}
\nuc{\wtn}{1,\cdots,n}
\nuc{\wtT}{1,\cdots,T}
\nuc{\woN}{\frac{1}{N}}

\nuc{\xt}{x_t}
\nuc{\Xt}{X_t}
\nuc{\Xsubi}{X_i}
\nuc{\Xsi}{X_i}

\nuc{\yt}{y_t}

\begin{document}

\title{State Space Methods for Granger-Geweke Causality Measures
\thanks{This work was partially supported by the NIH}}

\author{Victor Solo\thanks{v.solo@unsw.edu.au}\\    
School of Electrical Engineering\\
University of New South Wales, Sydney, AUSTRALIA}  

\date{January 16, 2015}                           
\maketitle

\begin{abstract}                          
At least two recent developments have put the
spotlight on some significant gaps in the theory of 
multivariate time series. 
The recent interest in the dynamics of networks;
and the advent, across a range of applications,
of measuring modalities
that operate on different temporal scales.

Fundamental to the description of network dynamics
is the direction of interaction between nodes, accompanied
by a measure of the strength of such interactions.
Granger causality (GC) and its associated frequency domain
strength measures (GEMs) (due to Geweke) provide
a framework for the formulation and analysis of these issues.
In pursuing this setup, three significant unresolved issues emerge.

Firstly computing GEMs 
involves computing submodels of
vector time series models, for which reliable methods do not exist;
Secondly the impact of filtering on GEMs has never been definitively
established. Thirdly the impact of downsampling on GEMs has never
been established. In this work, using state space methods,
we resolve all these issues
and illustrate the results with some simulations.
Our discussion is motivated by some problems
in (fMRI) brain imaging but is of general applicability.
\end{abstract}


%
%
\nuc{\anb}[2]{\ens{(#1 \mbox{ } #2)}}
\nuc{\abcd}[4]{\ens{\big({#1 \atop #2} \mbox{ } {#3 \atop #4}\big)}}
%
\nuc{\expj}[1]{exp(#1)}  
%
\nuc{\Lm}{L^{-1}} 
\nuc{\Lp}{L}     
\nuc{\Ab}{\bar{A}}
\nuc{\Abm}{\bar{A}_m}
\nuc{\As}{A_s}
\nuc{\at}{a_t}
\nuc{\bsubt}{b_t}
\nuc{\aum}{A^m}
\nuc{\cx}{C_X}
\nuc{\cy}{C_Y}
\nuc{\chisq}{\chi^2}
\nuc{\bx}{B_X}
\nuc{\by}{B_Y}
\nuc{\bo}{B_o}
%
%
\nuc{\ddash}{$"$}
\nuc{\dlamtpi}{\frac{d\lam}{2\pi}}

\nuc{\evec}{eigenvector }

\nuc{\bxt}{\bar{X}_t}
\nuc{\byt}{\bar{Y}_t}
\nuc{\bepxt}{\bar{\ep}_{X,t}}
\nuc{\bnt}{\bar{n}_t}
\nuc{\bzt}{\bar{Z}_t}

\nuc{\dlamotpi}{\frac{d\lam}{2\pi}}
\nuc{\calc}{{\mathcal C}}
\nuc{\calo}{{\mathcal O}}
\nuc{\dy}{d_y}
\nuc{\dx}{d_x}
\nuc{\dxx}{D_{XX}}
\nuc{\dxy}{D_{XY}}
\nuc{\dyx}{D_{YX}}
\nuc{\dyy}{D_{YY}}

\nuc{\epa}{\ep_a}
\nuc{\epb}{\ep_b}
\nuc{\est}{e_t}
\nuc{\epxt}{\ep_{X,t}}
\nuc{\epyt}{\ep_{Y,t}}
\nuc{\epbk}{\bar{\ep}_k}
\nuc{\epbxk}{\bar{\ep}_{X,k}}
\nuc{\epbyk}{\bar{\ep}_{Y,k}}
\nuc{\epxoyt}{({_{\ep_{X,t}} \atop ^{\ep_{Y,t}}})}
\nuc{\elam}{\expj{j\lam}}
\nuc{\elamm}{\expj{-j\lam}}

\nuc{\fyx}{F_{Y\ra X}}
\nuc{\fxy}{F_{X\ra Y}}
\nuc{\fyxh}{\hat{F}_{Y\ra X}}
\nuc{\fxyh}{\hat{F}_{X\ra Y}}
\nuc{\fbyx}{F_{\bar{Y}\ra \bar{X}}}
\nuc{\fbxy}{F_{\bar{X}\ra \bar{Y}}}
\nuc{\fyxlam}{f_{Y\ra X}(\lam)}
\nuc{\fxylam}{f_{X\ra Y}(\lam)}
\nuc{\fbyxlam}{f_{\bar{Y}\ra\bar{X}}(\lam)}
\nuc{\fbxylam}{f_{\bar{X}\ra\bar{Y}}(\lam)}
\nuc{\felam}{f_e(\lam)}
\nuc{\fxlam}{f_X(\lam)}
\nuc{\fylam}{f_Y(\lam)}
\nuc{\fbxlam}{f_{\bar{X}}(\lam)}
\nuc{\fbzlam}{f_{\bar{Z}}(\lam)}
\nuc{\fydx}{F_{Y.X}}
\nuc{\fxoy}{F_{XoY}}
\nuc{\fydxh}{\hat{F}_{Y.X}}
\nuc{\fxoyh}{\hat{F}_{XoY}}
\nuc{\fg}{f_G}
\nuc{\fy}{f_Y}
\nuc{\fxz}{f_X(\Lp)}
\nuc{\fez}{f_e(\Lp)}

\nuc{\fa}{f_a}
\nuc{\fb}{f_b}
%
\nuc{\Gc}{G_c}
\nuc{\Go}{G_o}
\nuc{\Goc}{G_{o,c}}

\nuc{\hxik}{\hat{\xi}_k}
\nuc{\hxikp}{\hat{\xi}_{k+1}}
\nuc{\hxx}{H_{XX}}
\nuc{\hxy}{H_{XY}}
\nuc{\hyx}{H_{YX}}
\nuc{\hyy}{H_{YY}}
\nuc{\he}{H_e}
\nuc{\hex}{H_{eX}}
\nuc{\hexz}{H_{eX}(\Lp)}
\nuc{\hbex}{H_{e\bar{X}}}
\nuc{\Hb}{\bar{H}}
\nuc{\gss}{(A,C,[Q,S,R])}
\nuc{\iss}{(A,C,B,\Sg_\ep)}
\nuc{\gamx}{\gam_x}
\nuc{\gamy}{\gam_y}
\nuc{\gzm}{g(z^{-1})}
\nuc{\Gzm}{G(z^{-1})}
%
\nuc{\heyxz}{H_{e,YX}(\Lp)}  
\nuc{\ho}{H_o}
\nuc{\hoex}{H_{oeX}}
\nuc{\hnt}{\hat{n}_t}
\nuc{\hyet}{\hat{Y}_{E,t}}
\nuc{\gamo}{\gam_0}
\nuc{\gamw}{\gam_1}
%
%
%
\nuc{\hxyz}{H_{XY}(L)}
\nuc{\hxxz}{H_{XX}(L)}
\nuc{\hyxz}{H_{YX}(L)}
\nuc{\hyyz}{H_{YY}(L)}

\nuc{\hc}{H_c}

\nuc{\hbxyz}{\bar{H}_{XY}(L)}
\nuc{\hbxxz}{\bar{H}_{XX}(L)}
\nuc{\hbyxz}{\bar{H}_{YX}(L)}
\nuc{\hbyyz}{\bar{H}_{YY}(L)}

\nuc{\kbex}{K_{e\bar{X}}}
\nuc{\km}{K_m^\ast}
\nuc{\Kbm}{\bar{K}_m}
\nuc{\lm}{L_m}
\nuc{\kcx}{K_{(X)}}
\nuc{\kcy}{K_{(Y)}}
\nuc{\Kt}{K_t}
\nuc{\Ks}{K_s}

\nuc{\Lum}{L\upm}
\nuc{\Jo}{J_o}
\nuc{\Jb}{\bar{J}}

\nuc{\ka}{k_a}
\nuc{\kb}{k_b}
\nuc{\kc}{k_c}

\nuc{\nuxk}{\nu_{X,k}}
\nuc{\nuyk}{\nu_{Y,k}}
\nuc{\nuxl}{\nu_{X,l}}
\nuc{\ozk}{\ol{z}_k}
\nuc{\oxik}{\ol{\xi}_k}
\nuc{\oxikp}{\ol{\xi}_{k+1}}
\nuc{\oepk}{\ol{\ep}_k}
\nuc{\nuk}{\nu_k}
\nuc{\nt}{n_t}
\nuc{\omx}{\Om_X}
\nuc{\phizm}{\phi(\Lp)}
\nuc{\Phizm}{\Phi(\Lp)}
\nuc{\pcx}{P_{(X)}}
\nuc{\phix}{\phi_x}
\nuc{\Phix}{\Phi_X}
\nuc{\phiy}{\phi_y}
\nuc{\phia}{\phi_a}
\nuc{\phib}{\phi_b}
\nuc{\Phiz}{\Phi(\Lp)}
\nuc{\phixyz}{\Phi_{XY}(\Lp)}
\nuc{\phixxz}{\Phi_{XX}(\Lp)}
\nuc{\phixx}{\Phi_{XX}}
\nuc{\phiyxz}{\Phi_{YX}(\Lp)}
\nuc{\phiyyz}{\Phi_{YY}(\Lp)}
\nuc{\nxt}{n_{X,t}}
\nuc{\nyt}{n_{Y,t}}

\nuc{\py}{p_Y}
\nuc{\Pt}{P_t}
\nuc{\Ptp}{P_{t+1}}
\nuc{\Pm}{P_m^\ast}

\nuc{\Phib}{\bar{\Phi}}
\nuc{\Psib}{\bar{\Psi}}

\nuc{\Sgep}{\Sg_{\ep}}
\nuc{\sgo}{\Sg_o}
\nuc{\sgoo}{\Sg^o}
\nuc{\sgep}{\Sg_{\ep}}
\nuc{\sgepx}{\Sg_{X,\ep}}
\nuc{\sgepy}{\Sg_{Y,\ep}}
\nuc{\sgepyx}{\Sg_{YX,\ep}}
\nuc{\sgepxy}{\Sg_{XY,\ep}}
\nuc{\sgepygx}{\Sg_{(Y|X),\ep}}
\nuc{\sgepxgy}{\Sg_{(X|Y),\ep}}
\nuc{\sgxe}{\sgepx}
\nuc{\sgyxe}{\sgepyx}
\nuc{\sgxye}{\sgepxy}
\nuc{\sgye}{\sgepy}
\nuc{\sgxem}{\Sg_{X,\ep}^{-1}}
\nuc{\qssr}{({_Q \atop ^{S^T}} {_{S} \atop ^R})}
\nuc{\qsm}{Q_m}
\nuc{\Qb}{\bar{Q}}
\nuc{\Qbm}{\bar{Q}_m}
\nuc{\Qs}{Q_s}
\nuc{\Qm}{Q_m}
\nuc{\Qmm}{Q_{m-1}}
\nuc{\thrd}{\frac{1}{3}}
\nuc{\ry}{r_Y}
\nuc{\rg}{r_G}
\nuc{\thty}{\tht_Y}
\nuc{\thtg}{\tht_G}
\nuc{\thtx}{\tht_x}
\nuc{\sga}{\sg_a}
\nuc{\sgb}{\sg_b}
\nuc{\sm}{S_m}
\nuc{\Sgb}{\overline{\Sg}}
%
\nuc{\Rb}{\overline{R}}
\nuc{\Sb}{\overline{S}}
\nuc{\stsp}{state space }
\nuc{\taua}{\tau_a}
\nuc{\taub}{\tau_b}

\nuc{\vm}{V_m^\ast}
\nuc{\vbm}{\bar{V}_m}
\nuc{\vcx}{V_{(X)}}
\nuc{\supa}{^a}

\nuc{\hxit}{\hat{\xi}_t}
\nuc{\hxitp}{\hat{\xi}_{t+1}}
\nuc{\wk}{w_k}
\nuc{\wt}{w_t}
\nuc{\vnu}{V_\nu}
\nuc{\Vt}{V_t}
\nuc{\Vtm}{V_t^{-1}}
\nuc{\upcm}{^{(m)}}
\nuc{\zimai}{(zI-A)^{-1}}
\nuc{\wophixzm}{\frac{1}{1-\phi_X(z^{-1})}}
\nuc{\wophiyzm}{\frac{1}{1-\phi_Y(z^{-1})}}
\nuc{\wxk}{w_{X,k}}
\nuc{\wyk}{w_{Y,k}}
\nuc{\wtxk}{\tilde{w}_{X,k}}
\nuc{\uph}{^{\half}}
\nuc{\xit}{\xi_t}
\nuc{\xitp}{\xi_{t+1}}
\nuc{\zmt}{z^{-2}}
\nuc{\zt}{z_t}
\nuc{\xoyt}{({_{x_t} \atop ^{y_t}})}
\nuc{\ybk}{\bar{y}_k}
\nuc{\xbk}{\bar{x}_k}
\nuc{\zbk}{\bar{z}_k}
\nuc{\xbat}{\bar{x}_t}
\nuc{\ybat}{\bar{y}_t}
\nuc{\zbt}{\bar{z}_t}
\nuc{\xabm}{X_{a,b}^{-}}
\nuc{\xabp}{X_{a,b}^{+}}
\nuc{\xoa}{X_a^0}
\nuc{\xab}{X_a^b}
\nuc{\xmt}{X_{-\infty}^t}
\nuc{\ymt}{Y_{-\infty}^t}
\nuc{\Xtp}{X_{t+1}}
\nuc{\xupt}{X^t}
\nuc{\xuptp}{X^{t+1}}
\nuc{\yupt}{Y^t}
\nuc{\yuptp}{Y^{t+1}}
\nuc{\xuptpp}{X^{t+p}}
\nuc{\yuptpp}{Y^{t+p}}
\nuc{\xuptppp}{X^{t+p+1}}
\nuc{\xtpupinf}{X_{t+1}^{\infty}}
\nuc{\xmtmm}{X^{-1}_{-(t-1)}}
\nuc{\ymtmm}{Y^{-1}_{-(t-1)}}
\nuc{\Ytp}{Y_{t+1}}
%
\nuc{\xb}{\bar{X}}
\nuc{\xix}{\xi_x}
\nuc{\xiy}{\xi_y}
\nuc{\zetak}{\zeta_k}
\nuc{\zimam}{(L^{-1}I-A)^{-1}}
\nuc{\zimaumm}{(L^{-1}I-A^m)^{-1}}
\nuc{\liima}{(L^{-1}I-A)^{-1}}
\nuc{\imal}{(I-AL)^{-1}}
\nuc{\liimam}{(L^{-1}I-A^m)^{-1}}
\nuc{\imaml}{(I-A^mL)^{-1}}
\nuc{\czl}{C_Z(L)}
\nuc{\czk}{C_{Z,k}}
\nuc{\czo}{C_{Z,0}}
\nuc{\alp}{\al_\perp}
\nuc{\cxl}{C_X(L)}
\nuc{\cyl}{C_Y(L)}
\nuc{\cyxl}{C_{YX}(L)}
\nuc{\bcyxl}{\bar{C}_{YX}(L)}
\nuc{\bczl}{\bar{C}_Z(L)}
\nuc{\bcxl}{\bar{C}_X(L)}
\nuc{\bepxs}{\bar{\ep}_{X,s}}
\nuc{\bcx}{\bar{C}_X}
\nuc{\hepxt}{\hat{\ep}_{X,t}}
\nuc{\hepxms}{\hat{\ep}_{X,ms}}

\nuc{\eptmk}{\ep_{t-k}}
\nuc{\dombtpi}{\frac{d\bar{\om}}{2\pi}}
\nuc{\czom}{C_Z(e^{j\om})}

\nuc{\epxms}{\ep_{X,ms}}
\nuc{\ejkom}{e^{jk\om}}
\nuc{\byes}{\bar{Y}_{E,s}}
\nuc{\bns}{\bar{n}_s}
\nuc{\bep}{\bar{\ep}}
\nuc{\hep}{\hat{\ep}}
\nuc{\bac}{\bar{C}}
\nuc{\bys}{\bar{Y}_s}
\nuc{\bxs}{\bar{X}_s}
\nuc{\bzs}{\bar{Z}_s}
\nuc{\byis}{\bar{Y}_{I,s}}

\nuc{\hx}{H(X)}
\nuc{\hxgy}{H(X|Y)}
\nuc{\hygx}{H(Y|X)}
\nuc{\ixy}{I(X;Y)}
\nuc{\hy}{H(Y)}
\nuc{\gz}{|Z}
\nuc{\xcy}{X;Y}
\nuc{\fzom}{f_Z(\om)}
\nuc{\gamzk}{\Gam_{Z,k}}
\nuc{\gamzkm}{\Gam_{Z,km}}
\nuc{\intpp}{\int_{-\pi}^\pi}
\nuc{\intotpi}{\int_0^{2\pi}}
\nuc{\fz}{f_Z}
\nuc{\fbz}{f_{\bar{Z}}}
\nuc{\fwom}{F_W(\omega)}
\nuc{\grcs}{{_{GC} \atop ^{\lra}}}
\nuc{\lgrcs}{{_{GC} \atop ^{\lra}}}
\nuc{\wgrcs}{{_{WGC} \atop ^{\lra}}}
\nuc{\wlgrcs}{{_{WGC} \atop ^{\lra}}}
\nuc{\wgrcsa}{{_{WGC^\ast} \atop ^{\lra}}}
\nuc{\sgrcs}{{_{SGC} \atop ^{\lra}}}
\nuc{\sgrcsa}{{_{SGC^\ast} \atop ^{\lra}}}
\nuc{\slgrcs}{{_{SGC} \atop ^{\lra}}}
\nuc{\usgrcs}{{_{USGC} \atop ^{\lra}}}
\nuc{\ugrcs}{{_{UGC} \atop ^{\lra}}}
\nuc{\ulgrcs}{{_{UGC} \atop ^{\lra}}}
\nuc{\haw}{\hat{\omega}}
\nuc{\hyes}{\hat{Y}_{E,s}}

\nuc{\lt}{$L_t$ }
\nuc{\lot}{$L_t^0$ }
\nuc{\omba}{\bar{\omega}}
\nuc{\omb}{\bar{\omega}}
\nuc{\limminf}{{_{lim} \atop ^{m\ra\infty}}}

\nuc{\rt}{$R_t$ }
\nuc{\rot}{$R_t^0$ }
\nuc{\sgx}{\Sg_X}
\nuc{\sgxy}{\Sg_{XY}}
\nuc{\sgyx}{\Sg_{YX}}
\nuc{\sgy}{\Sg_Y}
%

\nuc{\xw}{X_1}
\nuc{\xtwo}{X_2}
\nuc{\xtre}{X_3}
%
\nuc{\xat}{X_{\al,t}}
\nuc{\xbt}{X_{\be,t}}
\nuc{\ztt}{Z_{\tht,t}}
\nuc{\zbs}{\bar{Z}_s}
\nuc{\zms}{Z_{ms}}
\nuc{\xbs}{\bar{X}_s}
\nuc{\xms}{X_{ms}}
\nuc{\ybs}{\bar{Y}_s}
\nuc{\yms}{Y_{ms}}
\nuc{\wom}{\frac{1}{m}}
%
\nuc{\wscs}{{_{WSC} \atop ^{\lra}}}
\nuc{\sscs}{{_{SSC} \atop ^{\lra}}}
\nuc{\wscsa}{{_{WSC^\ast} \atop ^{\lra}}}
\nuc{\sscsa}{{_{SSC^\ast} \atop ^{\lra}}}
\nuc{\wpcs}{{_{WPC} \atop ^{\lra}}}
\nuc{\spcs}{{_{SPC} \atop ^{\lra}}}
\nuc{\wpcsa}{{_{WPC^\ast} \atop ^{\lra}}}
\nuc{\spcsa}{{_{SPC^\ast} \atop ^{\lra}}}
\nuc{\yet}{Y_{E,t}}
\nuc{\yit}{Y_{I,t}}
\nuc{\upw}{^{-1}}
\nuc{\yonm}{Y_1^{n-1}}
\nuc{\zdi}{Z_i}
\nuc{\zimw}{Z^{i-1}_1}
\nuc{\zmimw}{Z^{-1}_{-(i-1)}}
\nuc{\zwn}{Z_1^n}
\nuc{\xmimw}{X^{-1}_{-(i-1)}}
\nuc{\ymimw}{Y^{-1}_{-(i-1)}}
\nuc{\zdm}{Z^{-}}
\nuc{\xmmn}{X_{-n}^{-1}}
\nuc{\ymmn}{Y_{-n}^{-1}}
%
%

\nuc{\px}{p_x}
\nuc{\xon}{X_1^n}
\nuc{\xton}{X_1,\cdots,X_n}
\nuc{\xtn}{X_1,\cdots,X_n}
\nuc{\yon}{Y_1^n}
\nuc{\xonm}{X_1^{n-1}}
\nuc{\yonp}{Y_1^{n+1}}
\nuc{\yton}{Y_1,\cdots,Y_n}
\nuc{\zon}{Z_1^n}
\nuc{\fon}{\frac{1}{n}}
%
\nuc{\paperp}{\perp}
\nuc{\calm}{{\cal M}}
%
\nuc{\xupn}{X^n}
\nuc{\xupnp}{X^{n+1}}
\nuc{\xupinf}{X^{\infty}}
\nuc{\yupn}{Y^n}
\nuc{\yupnp}{Y^{n+1}}
\nuc{\yupinf}{Y^{\infty}}
\nuc{\xupnpp}{X^{n+p}}
\nuc{\yupnpp}{Y^{n+p}}
\nuc{\xupnppp}{X^{n+p+1}}
\nuc{\xnpupinf}{X_{n+1}^{\infty}}
\nuc{\xmnmm}{X^{-1}_{-(n-1)}}
\nuc{\ymnmm}{Y^{-1}_{-(n-1)}}

\nuc{\xo}{X^0}
\nuc{\yo}{Y^0}
\nuc{\zo}{Z^0}
\nuc{\xm}{X^{-}}
\nuc{\ym}{Y^{-}}
\nuc{\zum}{Z^{-}}
\nuc{\xp}{X^{+}}
\nuc{\yp}{Y^{+}}
\nuc{\zp}{Z^{+}}
\nuc{\comdots}{,\cdots,}
\nuc{\mn}{m_n}
\nuc{\xn}{x_n}
\nuc{\xnp}{x_{n+1}}
\nuc{\yn}{y_n}
\nuc{\ynp}{y_{n+1}}

\nuc{\ns}{ance}
\nuc{\ect}{empirical causality testing }
\nuc{\ectp}{empirical causality testing procedure }
\nuc{\evals}{eigenvalues }

\nuc{\dt}{discrete time }
\nuc{\ct}{continuous time }
\nuc{\biva}{bivariate }
\nuc{\casy}{causality }
\nuc{\anly}{analysis }
\nuc{\catch}{catchment }
\nuc{\cata}{catchment area }
\nuc{\cind}{conditional independence }
\nuc{\cent}{conditional entropy }
\nuc{\chr}{chain rule }
\nuc{\cdl}{conditional }
\nuc{\dcri}{{\bf DCRI} }
\nuc{\dri}{{\bf DRI} }
\nuc{\dpi}{{\bf DPI} }
\nuc{\cdip}{conditional independence }
\nuc{\Cdip}{Conditional independence }
\nuc{\dngc}{does not Granger cause }
\nuc{\dnsc}{does not Sims cause }
\nuc{\dn}{does not }

\nuc{\conl}{controllable }
\nuc{\conly}{controllability }
\nuc{\eval}{eigenvalue }


\nuc{\gc}{Granger causality }
\nuc{\grc}{Granger cause }
\nuc{\gct}{Granger causality testing }
\nuc{\gcs}{Granger causes }
\nuc{\gctp}{Granger causality testing procedure }
\nuc{\inft}{information theory }

\nuc{\mult}{multivariate }
\nuc{\nst}{nonstationary }
\nuc{\nsty}{nonstationarity }
\nuc{\obint}{observation interval }
\nuc{\mui}{mutual information }
\nuc{\Mui}{Mutual Information }
\nuc{\nl}{nonlinear }
\nuc{\nons}{nonstationary }
\nuc{\pnd}{purely non-deterministic }
\nuc{\obse}{observable }
\nuc{\obsy}{observability }
\nuc{\lamw}{\lam_1}
\nuc{\lamt}{\lam_2}

\nuc{\procs}{processes }
\nuc{\proc}{process }
\nuc{\stty}{stationarity }
\nuc{\sint}{sampling interval }
\nuc{\rivl}{river level }
\nuc{\pev}{prediction error variance }
\nuc{\pca}{principal component analysis }
\nuc{\pc}{principal component }
\nuc{\pdf}{probability density function }
\nuc{\pdfs}{probability density functions }
\nuc{\sym}{{\bf SYM} }
\nuc{\scas}{Sims causality }

\nuc{\wod}{Wold decomposition }
\nuc{\varie}{variance }

\section{\bf Introduction}
Following the operational development of the notion
of causality by \cc{GRAN69} and \cc{SIMS72},
Granger causality (henceforth denoted GC) 
analysis has become an important
part of time series and econometric testing and inference
e.g. \cc{HAML94}. It has also been applied in the biosciences,
\cc{KAMI91}, \cc{BERN99}, \cc{DING00b};
climatology
(global warming) \cc{SUNA96}, \cc{STER97}, \cc{TRIA05};
and most recently
functional magnetic resonance imaging (fMRI).

Since its introduction into fMRI 
\cc{GOEB05}, \cc{OZAK05} it has become the subject of
an intense debate: e.g. see \cc{ROEB11} and 
associated commentary
on that paper. 
There are two main issues in that debate but
which occur more widely in dynamic networks.
Firstly,
the impact of downsampling on GC.
In the fMRI neuro-imaging application
causal processes
may operate on a time-scale of order 
tens of milli-seconds
whereas the recorded
signals are only available on a one-second time-scale.
So it is natural to wonder if GC analysis on a slow
time-scale can reveal dynamics on a much faster time-scale.
Secondly,  the impact of filtering on GC
due to the
hemodynamic response function
which relates the neural activity to the recorded fMRI
signal.
Since intuitively GC will be sensitive to
time delay, 
the variability of the hemodynamic response
function, particularly spatially varying 
time to onset and time to peak (confusingly
called delay in the fMRI literature)
has been suggested as a potential
source of problems \cc{DESP10a},\cc{HANW12}.

An important advance in GC theory and tools was made by \cc{GEWK82}
who provided measures of the strength of causality (henceforth
called GEM for Geweke causality measure)
including frequency domain decompositions of them. 
Subsequently
it was pointed out that the GEMs are measures of mutual information
\cc{WAXR87}.
The GEMs were extended to
conditional causality in \cc{GEWK84}. 
However GEMs have not found as wide application
as they should have, partly because of some technical difficulties
in calculating them discussed further below. But GEMs 
(and their frequency domain versions) are
precisely the tool needed to pursue both the 
GC downsampling and filtering questions.


In the econometric literature,
it was appreciated early that downsampling, especially
in the presence of aggregation could cause problems.
This was implicit in work of \cc{SIMS71}
; mentioned also in
work of \cc{CREI87} who gave an example
of contradictory causal analysis based
on monthly versus quaterly data and also
discussed in \cc{MARC85}. But precise general conditions
under which problems do and do not arise have
never been given. We do so below.

Some of the above econometric discussion is framed in terms
of sampling of continuous time models \cc{SIMS71},
\cc{MARC85},\cc{CREI87}. And authors such as \cc{SIMS71}
have suggested that models are best formulated initially
in continuous time. While this is a view the author
has long shared we deal with only discrete time models
here. To cast our development in terms of continuous
time models would require a considerable development
of its own without changing our basic message.

The issue at stake, in its simplest form, is the following. 
Suppose that a pair of (possibly vector) processes posses a 
unidirectional GC relation but suppose
measurements are only available at a slower time-scale
on filtered series.
Then two questions arise.
The first, which we call the \ul{forward} question, is this:
Is the unidirectional Granger causal relation preserved? 
The second, which we call the \ul{reverse} question, is
harder. Suppose the downsampled filtered
series exhibit a uni-directional
GC relation; does that mean the underlying 
unfiltered
faster 
time-scale processes do?
The latter question is the more important and so far
has received no theoretical attention.


In order to resolve these issues we need to develop
some theory and some computational/modeling tools. 
Firstly to compute GEMs
one needs to be able to find submodels 
from a larger 
(i.e. one having more time series) model. 
Thus to compute the GEMs between time series $\xt,\yt$
\cc{GEWK82},\cc{GEWK84} attempted to avoid this by fitting
submodels separately to $\xt$ to $\yt$ and 
then also fitting a joint model to $\xt,\yt$.
Unfortunately this can generate negative values for some of the
frequency domain GEMs \cc{DING06}. 
Properly computing submodels will
resolve this problem and previous work 
has not accomplished this
(we discuss the attempts in \cc{DUFT10} and \cc{DING06} below).

Secondly one needs to be able to compute how models
transform when downsampled. 
This has only been done in special cases 
\cc{PANW83}
or by methods
that are not computationally realistic. 
We provide
computationally reliable, state space based methods 
for doing this here.

Thirdly we need to study the effect of filtering on GEMs.
And then using these tools
one can compute filtered downsampled GEMs and hence
study the effect of sampling and filtering on GEMs.

To sum up we can say that
previous discussions including those above as well
as \cc{GEWK78},\cc{TELS67},\cc{SIMS71},\cc{MARC85}
fail to provide general algorithms for finding
submodels or models induced
by downsampling.
Indeed both these problems have remained
open problems in multivariate time series in their own right
for several decades and we resolve them here.
Further there does not seem to have been any 
theoretical discussion
of the effect of filtering on GEMs and we 
resolve that here also. To do that it turns
out that state space models provide the proper framework.

Throughout this work we deal with the dynamic interaction between
two vector time series. It is well known 
that if there is
a third vector time series involved in the dynamics but not accounted
for then spurious causality can occur
for reasons that have nothing to do with downsampling.
This situation has been
discussed by \cc{HSIAO82}; see also \cc{GEWK84}. Other causes of spurious
causality such as observation noise are also not discussed.
Of course the impact of downsampling in the presence of a third (vector)
variable is also of interest but will be pursued elsewhere.

Finally our whole discussion is carried out in
the framework of linear time series models.
It is of great interest to pursue nonlinear versions
of these issues but that will be a major task.

The remainder of the paper is organized as follows.
In section 2 we review and modify some 
state space results important for
system identification or model 
fitting and needed in the following sections.
In section 3 we develop \stsp methods
methods for computing submodels of innovations state space models.
In section 4 we develop
methods for transforming state space models under downsampling.
In section 5 we review 
GC and
GEMs and extend them to a state space setting.
In section 6 we study the effect of filtering on GC via frequency
domain GEMs.
In section 7 we give theory
to explain when causality is preserved under
downsampling. 
In section 8 
we discuss the reverse problem 
showing how spurious causality can be induced by
downsampling.
%
%
Conclusions are in section 9. There are three appendices.

\subsection{\bf Acronyms and Notation}
GC is Granger causality or Granger causes.
We use the GC designator alone where we make statements of 
interest in both weak and strong cases.
dn-gc is does not Granger cause;
GEM is Geweke causality measure;
SS is state space or state space model;
ISS is innovations state space model;
VAR is vector autoregression;
VARMA is vector autoregressive moving average process;
wp1 is with probability 1.

$\xab$ denotes the values $x_a,x_{a+1},\cdots,x_b$; so $X_a^a\equiv x_a$.
For stationary processes we have $a=-\infty$.
$\zm=L$ is the lag or backshift operator;
LHS denotes left hand side etc.
If $M,N$ are positive semi-definite matrices
then $M\geq N$ means $M-N$ is positive semi-definite.
A square matrix is stable if all its \evals have modulus $<1$.

 \section{\bf State Space}
\sco
The computational methods
we develop rely on state space techniques
and spectral factorization.

There is an intimate relation between
the steady state Kalman filter and
spectral factorization
which is fundamental to
our computational procedures.

So
in this section we review 
and modify
some basic results in 
state space theory, Kalman filtering
and spectral factorization.

In the sequel we deal with
two vector time series, which we collect together as,
$\zt=(\xt^T, \yt^T)^T$.

\subsection{\bf State Space Models}
We consider a general constant parameter SS model,
\EQN
\xitp=A\xit+w_t ,\;  \zt=C\xit+v_t                 \lab{ss}
\ENN
with positive semi-definite noise covariance,
$var\aob{w_t}{v_t}=\aobcod{Q}{S^T}{S}{R}$.
We refer to this as a SS model with parameters
(A,C,[Q,R,S]). 

It is common with SS models to take $S=0$, but for
equivalence between the class of VARMA models 
and the class of state space models it is necessary 
to allow $S\neq 0$.

Now by matrix partitioning,
$|\aobcod{Q}{S^T}{S}{R}|=|R||\Qs|,\Qs=Q-SR\upm S^T$.
So introduce,

\noi{\bf Noise Condition N}. 
$R$ is positive definite.

whereupon $\Qs$ is positive semi-definite.

\subsection{\bf Steady State Kalman Filter,
Innovations State Space (ISS) Models and 
the Discrete Algebraic Ricatti Equation (DARE)}
We now recall the Kalman filter 
for mean square estimation of the unobserved state sequence
$\xit$ from the observed time series $\zt$. It is given by
\cc{KAIL00}(Theorem 9.2.1),
\EQ
\hxitp=A\hxit+K_t\est,\;\;
\est=\zt-C\hxit,\mbox{ or } \zt=C\hxit+\est
\EN
where $\est$ is the
innovations sequence of variance $\Vt=R+C\Pt C^T$ and 
$K_t=(A\Pt C^T+S)\Vtm$ is the Kalman gain sequence 
and $\Pt$ is the state error variance matrix generated
from the Ricatti equation,
$\Ptp=A\Pt A^T+Q-\Kt\Vt \Kt^T$.

The Kalman filter is a time-varying filter
but we are interested in its steady state.
If there is a steady state 
i.e. $\Pt\ra P$ as $t\rai$ then
then the limiting state error
variance matrix $P$ will obey the so-called
discrete algebraic Ricatti equation ({\bf DARE})
\EQN
P=APA^T+Q-KVK^T,              \lab{dare}
\ENN
where $V=R+CPC^T$ and $K=(APC^T+S)V\upm$ is the corresponding
steady state Kalman gain. With some clever algebra
\cc{KAIL00}(section 9.5.1) the DARE can be rewritten
(the Ricatti equation can be similarly rewritten),
\EB
P=\As P\As^T+\Qs-\Ks V\Ks^T
\EE
where $\As=A-SR\upm C$ and $\Ks=\As PC^TV\upm$.

We now introduce two assumptions.\\
\noi{\bf Stabilizability Condition St}: 
$\As,\Qs\uph$ is stabilizable (see Appendix A)\\
\noi{\bf Detectability Condition De}: 
$\As,C$ is detectable.

\noi In Appendix A
it is shown this is equivalent to $A,C$ being detectable.
And also it holds automatically if $A$ is stable.\\

The resulting steady state Kalman filter can be written as,
\EQN
\hxitp=A\hxit+K\ept,\;\;        \zt=C\hxit+\ept            \lab{iss}
\ENN
where $\ept$ is the steady state innovation 
process and has variance
matrix $V$ and Kalman gain $K$.
This steady state filter provides a new state 
space representation
of the data sequence. We refer to it as an 
innovations
state space (ISS) model with parameters $(A,C,K,V)$. 
We summarize this in,

\noi{\bf Result I}. 
Given the SS model (\ref{ss}) with parameters $(A,C,[Q,R,S])$, 
then provided N,St,De hold:

(a)
The corresponding ISS model (\ref{iss}) 
with parameters $(A,C,K,V)$ 
can be found by solving the DARE 
(\ref{dare}) which has a unique
positive definite solution $P$.
%

(b)
$V$ is positive definite,
$(A,C)$ is detectable and
$A-KC$ is stable so that $(A,K)$ is controllable.
%

{\it Proof}. See appendix A.\\
%
{\it Remarks}.

(i) Henceforth an ISS model 
with parameters $(A,C,K,V)$ will be required to have
$V$ positive definite, $(A,C)$ detectable and 
$(A,K)$ controllable so that $A-KC$ is stable.

(ii) It is well known that any VARMA 
model can be represented as an
 ISS model and vice versa 
\cc{SOLO86},\cc{HDE88}.
%

(iii) Note that the ISS model with parameters $(A,C,K,V)$ can also be written
as the SS model with parameters $(A,C,[KVK^T,V,KV])$.

(iv) The DARE is a quadratic matrix equation but can be
computed using the
(numerically reliable) DARE command in matlab as follows.
Compute:
$[P,L_0,G]=DARE(A^T,C^T,Q,R,S,I)$
and then,
$V=R+CPC^T,K=G^T$.

(v) Note that stationarity is not required for this result.

\subsection{\bf Stationarity and Spectral Factorization}
Given an ISS model with parameters $(A,C,B,\Sg_\ep)$,
we now introduce,\\
\noi{\bf Condition Ev}: 
$A$ has all \evals with modulus $<1$ i.e. $A$ is a stability matrix.\\
With this assumption we can obtain an infinite vector moving average (VMA)
representation, an infinite vector autoregressive (VAR)
representation and a spectral factorization.
The following result is based on \cc{KAIL00}[Theorem 8.3.2] and surrounding
discussion.\\
\noi{\bf Result II}. For the ISS
model $(A,C,B,\Sg_\ep)$ obeying condition Ev we have,

\noi(a) Infinite VMA or Wold decomposition,
\EQN
\zt=H(L)\ept=(C\liima B+I)\ept =(C\imal B L+I)\ept        \lab{issb}
\ENN 
%
(b)
Infinite VAR representation,
\EQN
\ept=G(L)\zt=[I+C(L\upm I -A+KC)\upm K]\zt              \lab{ivar}
\ENN
%
(c) Spectral factorization. 
Put $L=\exp(-j\lam)$ then, $\zt$
has positive definite spectrum with spectral factorization as follows,
\EQN
f_Z(\lam)=[C\liima , I] 
\aobcodL{Q}{S^T}{S}{R} \abL{(L I-A^T)\upm C^T}{I}  \lab{spf}
=H(L)\sgep H^T(L\upm)
\ENN
(d) $H(L)$ is minimum phase i.e. its
inverse exists and is causal and stable.\\

\noi{\it Proof}. 
\ul{(a)}. Just write (\ref{iss}) in operator form.
The series is convergent wp1 and in mean square since $A$ is stable.

\ul{(b)}. Rewrite (\ref{iss}) as,
$\hxitp=(A-KC)\hxit+K\zt,\ept=\zt-C\hxit$. Then write this in operator
form. The series is convergent wp1 and in mean square since $A-KC$ is stable
and $\zt$ is stationary.

\ul{(c)}. 
Follows from standard formulae 
for spectra 
of filtered stationary time series applied to (a).

\ul{(d)}.
From (a),(b) $G(L)=H\upm(L)$
and by (b) $G(L)$ is causal and stable and the result follows.\\

\noi{\it Remarks}.

(i)
For further discussion of minimum phase 
filters see \cc{GREE88},\cc{SOLO86}.

(ii) Result II is a special case of a general result
that given a full rank multivariate spectrum $f_Z(\lam)$ there exists
a unique causal stable minimum phase spectral factor $H(L)$ 
with $H(0)=I$
and
positive definite innovations variance matrix
$\sgep$ such that (\ref{spf}) holds \cc{HDE88},\cc{GREE88}.
In general $det H(L)$ may have
some roots on the unit circle \cc{HANP88},\cc{GREE88}
but the assumptions in result II rule this case out. 
Such roots mean that some linear combinations of $\zt$ can be perfectly
predicted from the past \cc{HANP88},\cc{GREE88} something that is not
realistic in the fMRI application.

(iii) Result II is also crucial from a system identification
or model fitting point of view. From that point of view all we can
know (from second order statistics) is the spectrum and so
if, naturally, we want a unique model,
the only model we can obtain is the causal stable minimum phase model
i.e. the ISS model. The standard approach to SS model fitting
is the so-called state space subspace method \cc{DEIS95},\cc{BAUR05}
and indeed it delivers an ISS model. The alternative approach
of fitting a VARMA model \cc{HDE88},\cc{LUTK93} is equivalent to getting an ISS model.

(iv) We need result I however since when we form submodels
we do not immediately get an ISS model, rather we must compute it.
\section{\bf Submodels}
\sco
Our computation of causality measures requires
that we compute induced submodels.
In this section we show how to obtain a ISS
submodel from the ISS joint model.

Now we partition $\zt=(\xt,\yt)^T$ into subvector signals of interest
and partition the \stsp model correspondingly,
$C=\aob{\cx}{\cy}$ and $B=(\bx,\by)$.
%
We first read out a SS submodel for $\xt$ from the ISS model
for $\zt$. We have simply 
$\xitp=A\xit+\wt,\;\; \xt=\cx\xit+\epxt$
where,
$\wt=B\epxoyt=\bx\epxt+\by\epyt$.
We need to calculate the covariance matrix,
$var\aob{\wt}{\epxt}=\abcd{Q}{\Sb^T}{\Sb}{\Rb}$.
We find, $\Rb=\sgepx$,
$Q=var(\wt)=B\sgep B^T$ and,
$\Sb=E(\wt\epxt^T)=B\aob{\sgepx}{\sgepyx}=\bo$.
This leads to

\noi{\bf Theorem I}.
 Given the joint ISS model (\ref{iss}) or (\ref{issb}) for 
$\zt$, then under condition Ev,
the corresponding ISS submodel for $\xt$ namely
($A,\cx,\kcx,\omx)$ 
(the bracket notation $\kcx$ is used to avoid confusion
with e.g. $\cx,\sgepx$ which are submatrices)
can be found
by solving the DARE  (\ref{dare})
with $[Q,\Rb,\Sb]=[B\sgep B^T,\sgepx,B_o]$.

\noi{\it Proof}. 
Firstly we note by partitioning $|\sgep|=|\sgepx||\sgepygx|$
where $\sgepygx=\sgepy-\sgepyx\sgepx\upm\sgepxy$
so that $\sgepx$ and $\sgepygx$ are both positive definite.
Now we need only check conditions N,St,De of result I.
We need to show, 
$\Rb=\sgepx$ is positive definite,$(A,\cx)$ is detectable and
$(A-\Sb\Rb\upm \cx,(B\sgep B^T-\Sb\Rb\upm\Rb^T)\uph$ is 
stabilizable; in fact we show it is controllable.

The first is already established.
The second follows trivially since $A$ is stable.
We use the PBH test (see Appendix A) to check the third.

Suppose controllability fails, then by the PBH test, there
exists $q\neq 0$ with
 $\lam q^T=q^T(A-\Sb\Rb\upm \cx)$
and 
$0=q^T(B\sgep B^T-\Sb\Rb\upm\Sb^T)\uph
\Ra 0=q^T(B\sgep B^T-\Sb\Rb\upm\Sb^T)
=q^T(B\sgep B^T-B\aob{\sgepx}{\sgepyx}\sgepx\upm[\sgepx,\sgepxy]B^T)
=q^T(\bx,\by)\aobcod{0}{0}{0}{\sgepygx}\aob{\bx^T}{\by^T}\\
=q^T\by\sgepygx\by^T
\Ra 0=q^T\by\sgepygx\by^Tq
\Ra \pa\by^Tq\pa=0\Ra q^T\by=0$ 
since $\sgepygx$ is positive definite.
But then,
$\lam q^T=q^T(A-(\bx,\by)\aob{\sgepx}{\sgepyx}\sgepx\upm\cx)
=q^T(A-(\bx\cx+\by\sgepyx\sgepx\upm\cx)
=q^T(A-\bx\cx)$.
Thus $(A-\bx\cx,\by)$ is not controllable.
But this is a contradiction since we can find a matrix
,namely $\cy$ so that $A-\bx\cx-\by\cy=A-BC$ is stable.

\noi{\it Remarks}.

(i) For implementation in matlab
positive definiteness in constructing $Q$ can be an issue.
A simple resolution is to use a Cholesky factorization of 
$\sgep=L_\ep L_\ep^T$ and form $B_\ep=BL_\ep$
and then form $Q=B_\ep B_\ep^T$.

(ii) The $\pcx$ matrix from DARE,\\
$[\pcx,L_0,G]=DARE(A^T,C_X^T,B\sgep B^T,\sgepx,B_o,I)$
obeys
$\pcx=A\pcx A^T+B\sgep B^T-\kcx\omx\kcx^T$ and then
$\kcx=(A\pcx\cx^T+B_o)\omx\upm,\omx=\sgepx+\cx\pcx\cx^T$.

(iii) \cc{DUFT10} discuss a method for obtaining submodels
but it is flawed.  
Firstly it requires
the computation of the inverse of the
VAR operator. 
While this might be feasible (analytically) on a toy example,
there is no known numerically reliable
way to do this in general
(computation of determinants is notoriously ill-conditioned).
Secondly it requires the solution of simultaneous
quadratic autocovariance 
equations to determine VMA parameters for which no algorithm
is given. In fact these are precisely the equations required
for a spectral factorization of a VMA process.
There do exist reliable algorithms for doing this but given
the flaw already revealed 
we need not discuss this approach any further.\\

Next we state an important corollary:\\
\noi{\bf Corollary I}. 
Any submodel is in general a VARMA
model not a VAR. To put it another way the class of VARMA models
is closed under the forming of submodels whereas
the class of VAR models is not.

This means that VAR models are not generic and is a strong
argument against their use. Any vector time series can be
regarded as a submodel of a larger dimensional time series
and thus must in general obey a VARMA model. This result
(which is well known in time series folk lore)
is significant for econometrics where VAR models are in widespread use.

For the next section we need,\\
\noi{\bf Theorem II}.
For the joint ISS model (\ref{iss})  or (\ref{issb}) for $\zt$
with conditions St,De holding and
with induced submodel for $\xt$ given in Theorem I,
we have,
\EBN
\fxlam&=&h_X(\elamm)\omx h_X^T(\elam)                  \lab{fxlam}\\
h_X(L)&=&I+\cx\liima \kcx =I+L \cx\imal\kcx \nn\\
ln|\omx|&=&\int_{-\pi}^{\pi}ln|\fxlam|\dlamtpi\nn
\EEN

\section{\bf Downsampling}
\sco
There are two approaches to the problem of finding
the model obeyed by a downsampled process;
frequency domain and time domain. 
While the
the general formula
for the spectrum of a sampled process has long been known,
it is not straightforward to use
and has not yielded any general computational
approach to finding submodels of parameterized spectra.
Otherwise the most complete (time domain) work seems to be that of 
\cc{PANW83} who only treat the first and second
order scalar cases. There is work in the engineering
literature for systems with observed inputs but that
is also limited and in any case not helpful here.
We follow a SS route.

We begin with the ISS model (\ref{iss}).
Suppose we downsample the observed signal $\zt$ with 
sampling multiple $m$.
Let $t$ denote the fine time scale and $k$ the 
coarse time scale so
$t=mk$. The downsampled signal is $\ozk=z_{mt}$.
To develop the SS model for $\ozk$ we iterate 
the SS model above to
obtain
\EQ
\xi_{t+l}=A^l\xit+\Sg_1^l A^{l-i}B\ep_{t+i-1}
\EN
Now set $t=mk,l=m$ and denote sampled signals,
$\oxik=\xi_{mk},\ozk=z_{mk},\oepk=\ep_{mk}$.
Then we find,
\EQ
\oxikp=A^m\oxik+\wk\;\;,
\ozk=C\oxik+\oepk
\EN
where $\wk=\Sg_1^m A^{m-i}B\ep_{km+i-1}$.
We now use result I to find the ISS model 
corresponding to this
SS model.

We have first to calculate the model covariances,
\EBN
E(\oepk\oepk^T)&=&\sgep=R \nn\\
E(\wk\oepk^T)&=&A^{m-1}B\sgep=\sm                        \lab{isssa}\\
E(\wk\wk^T)&=&Q_m   \nn\\
Q_m&=&\Sg_1^m A^{m-i}B\sgep B^T(A^T)^{m-i} 
=\Sg_0^{m-1}A^r B\sgep B^T (A^T)^r          \lab{isssq}\\
\Ra Q_m&=&AQ_{m-1}A^T+B\sgep B^T, m\geq 2, Q_1=B\sgep B^T\nn
\EEN
We now obtain,

\noi{\bf Theorem III}.
Given the ISS model (\ref{iss}), then under condition E,
for $m>1$, the ISS
model for the downsampled process $\ozk=z_{mt}$ is
$(A^m,C,\km,\vm)$ obtained by solving the DARE
with SS model $(A^m,C,[Q_m,R,\sm]$) where $Q_m$ 
is given in (\ref{isssq})
and $\sm$ is given in (\ref{isssa}).

\noi{\it Proof}. Using result I we need to show the following.
$R$ is positive definite, $(A^m,C)$ is detectable
and $(A^m-\sm R\upm C),(Q_m-\sm R\upm \sm^T)\uph)$
is stabilizable; in fact we show controllability.
The first holds trivially; the second also since $A$
is stable and thus so is $A^m$. For the third
we use the PBH test.

Suppose \conly fails. Then there is a left \evec $q$
(possibly complex) 
with $\lam q^T=q^T(A^m-\sm R\upm C)=q^T A^{m-1}(A-BC)$ 
and 
$q^T(Q_m-\sm R\upm \sm^T)
=0=q^T \Sg_0^{m-2}A^r B\sgep B^T(A^T)^r
\Ra \Sg_0^{m-2} q^H \Sg_0^{m-1}A^r B\sgep B^T(A^T)^rq=0$.
Since $\sgep$ is positive definite this delivers
$\Sg_0^{m-2}\pa B^T(A^T)^rq\pa^2=0\Ra B^T(A^T)^rq=0$
for $r=0,\cdots,m-2$.

Using this, we now find,
$\lam q^T=q^T A^{m-1}(A-BC)
=q^T A^{m-2}(A-BC)^2+q^T A^{m-2}BC(A-BC)
=q^T A^{m-2}(A-BC)^2$.
Iterating this yields, $\lam q^T=q^T(A-BC)^m$.
Thus if $\lam_m$ is an $m$-th root of $\lam$
then $\lam_m q^T=q^T(A-BC)$. Since also $q^TB=0$
we thus conclude $(A-BC,B)$ is not controllable.
But this is a contradiction since,
$(A-BC)+BC=A$ is stable.

\noi{\it Remarks}. 

(i) In matlab we would compute,
$[\Pm,L_0,G_m]=\mbox{DARE}((A^m)^T,C^T,Q_m,\sgep,\sm,I)$,
yielding $\vm=\sgep+C\Pm C^T$ and $\km=G_m^T$.

(ii) More specifically $\Pm$ ($m>1$) obeys,
$\Pm=\Abm \Pm \Abm^T+\Qbm-\Kbm\vm\Kbm^T$,\\
where $\Abm=A^m-\sm\sgep\upm C=A^{m-1}(A-BC)$;
$\Qbm=\Qm-\sm\sgep\upm\sm^T=\Qmm;
\vm=\sgep+C\Pm C^T;\km=\Abm\Pm C^T(\vm)\upm$.


\section{\bf Granger Causality}
\sco
In this section we review 
and extend some basic results in 
Granger causality. In particular
we extend GEMs to the state space setting
and show how to compute them reliably.

Since the development of Granger causality
it has become clear \cc{DUFR98},\cc{DUFT10} 
that in general one cannot
address the causality issue with only one step ahead
measures as commonly used; one needs to look at causality
over all forecast horizons.
However one step measures are sufficient when one is considering
only two vector time series as we are \cc{DUFR98}[Proposition 2.3].

\subsection{\bf Granger Causality Definitions}
Our definitions 
of one step Granger causality
naturally draw on \cc{GRAN63},
\cc{GRAN69},\cc{SIMS72},\cc{SOLO86}
but are also influenced by
\cc{CAIN76b}, who, drawing on work of \cc{PIERH77}, distinguished
between weak and strong GC or what Caines calls weak and strong feedback
free processes.
We introduce:\\
%
\noi{\bf Condition WSS}: 
The vector time series $\xt,\yt$ are jointly second order 
stationary.

\noi{\bf Definition: Weak Granger Causality}.\\
 Under WSS, we say $\yt$ 
does not \ul{weakly} \grc (dn-wgc) $\xt$ if, for all $t$
\EQ
var(\Xtp|\xmt,\ymt)=var(\Xtp|\xmt)
\EN
Otherwise we say $\yt$ weakly \gcs (wgc) $\xt$.

Because of the elementary identity,
$var(X|Z)=E[var(X|Z,W)]+var[E(X|Z,W)]
=E[var(X|Z,W)]+E[(E(X|Z,W)-E(X|Z))(E(X|Z,W)-E(X|Z))^T|Z]$
the equality of variance matrices in the definition
also ensures the equality of predictions,
$E(\Xtp|\xmt,\ymt)=E(\Xtp|\xmt)$.

This definition agrees with \cc{GRAN69},\cc{CAIN75}
who do not use the designator weak and \cc{CAIN76b},\cc{SOLO86} who do.

\noi{\bf Definition: Strong Granger Causality}.\\
 Under WSS, we say $\yt$ 
does not \ul{strongly} \grc (dn-sgc) $\xt$ if, for all $t$,
\EQ
var(\Xtp|\xmt,\ymt,\Ytp)=var(\Xtp|\xmt)
\EN
Otherwise we say $\yt$ strongly \gcs (sgc) $\xt$.
Again equality of the variance matrices ensures equality
of predictions, $E(\Xtp|\xmt,\ymt,\Ytp)=E(\Xtp|\xmt)$.

This definition agrees with \cc{CAIN76b} and \cc{SOLO86}.

\noi{\bf Definition: FBI}. Feedback Interconnected.\\
If $\xt$ \gcs  $\yt$ and $\yt$ \gcs $\xt$ 
then we say $\xt,\yt$ are feedback interconnected.\\
%
\noi{\bf Definition: UGC}. Unidirectionally Granger Causes.\\ 
If $\xt$ \gcs $ \yt$ but $\yt$ dn-gc $\xt$ we say 
$\xt$ unidirectionally Granger causes $ \yt$.

\subsection{\bf Granger Causality for Stationary State Space Models}

Now we partition $\zt=(\xt,\yt)^T$ into subvector signals of interest
and partition the vector MA or \stsp model (\ref{issb})
correspondingly,
\EBN
%
%
%
%
%
\abL{{\xt}}{{\yt}}&=&                      
[I+
\abL{{\cx}}{{\cy}}
\zimam
\anb{{\bx}}{{\by}}
\abL{{\epxt}}{{\epyt}}                         \lab{issc}\\
%
&=&\aobcodL{{\hxxz}}{{\hyxz}}{{\hxyz}}{{\hyyz}}
\abL{{\epxt}}{{\epyt}}                           \lab{issd}\\
\Sg_\ep&=&var\abL{{\epxt}}{{\epyt}}
=\aobcodL{\sgxe}{{\sgyxe}}{{\sgxye}}{{\sgye}} \nn\\
\aobcodL{{\hxxz}}{{\hyxz}}{{\hxyz}}{{\hyyz}}
&=&\aobcodL{\cx\zimam\bx+I}{\cy\zimam\bx}{\cx\zimam\by}{\cy\zimam\by+I}\nn
\EEN
Now we recall results of \cc{CAIN76b}:\\
\noi{\bf Result III}: 
If $\zt=\aob{\xt}{\yt}$ obeys a Wold model of the form
$\zt=H_Z(L)\ept$ where $H_Z(L)$ is a one-sided square summable
moving average polynomial with $H_Z(0)=I$
which is partitioned as in (\ref{issd})
then:

(a) $\yt$ dn-wgc $\xt$ iff $\hxyz=0$.

(b) $\yt$ dn-sgc $\xt$ iff  $\hxyz=0$ and $\sgxye=0$. 

We can now state a new SS version of this result:\\
\noi{\bf Theorem IV}. 
For the stationary ISS model (\ref{issc},\ref{issd}):

(a) $\yt$ dn-wgc $\xt$ iff $\cx A^r \by=0, r\geq 0$. 

(b) $\yt$ dn-sgc $\xt$ iff $\cx A^r \by=0, r\geq 0$ and  $\sgxye=0$. 

\noi{\it Proof}. Follows immediately from result III since
$\hxyz=\Sg_0^\infty \cx A^r\by L^{r+1}$.\\

\noi{\it Remarks}.

(i) By the Cayley Hamilton Theorem we can replace (a)
with: $\cx A^r\by=0,0\leq r\leq n-1,n=dim(\xit)$.

(ii) Collecting these equations together gives
$\cx(\by,A\by,\cdots,A^{n-1}\by)=0$ which says that 
the pair $(A,\by)$ is not controllable. Also we have,\\
$\by^T(\cx^T,A^T\cx^T,\cdots,(A^T)^{n-1}\cx^T)=0$ which
says that the pair $(\cx,A)$ is not observable.
Thus the representation of $\hxyz$ is not minimal.

From a data analysis point of view we need
 to embed this result
in a well behaved hypothesis test. Results of \cc{GEWK82}, suitably modified,
allow us to do this.

\subsection{\bf Geweke Causality Measures for SS Models}

Although much of the discussion in \cc{GEWK82} 
is in terms of VARs we can show it applies more generally.
We begin as \cc{GEWK82} did with the following definitions.
Firstly,
$\fyx=ln\frac{|\Om_X|}{|\sgxe|}$ is a 
measure of the gain in using
the past of $y$ to predict $x$ beyond using just the past of $x$;
similarly introduce $\fxy=ln\frac{|\Om_Y|}{|\sgye|}$. Next define the instantaneous
influence measure, $\fydx=ln\frac{|\sgxe||\sgye|}{|\Sg_\ep|} $.
These are then joined in the
fundamental decomposition \cc{GEWK82},
\EQN
\fxoy=\fyx+\fxy+\fydx   \lab{decomp} 
\ENN
where, $\fxoy=ln\frac{|\Om_X||\Om_Y|}{|\sgep|}$.
%
\cc{GEWK82} then proceeds to decompose these measures in the frequency domain.
Thus the frequency domain GEM for the dynamic influence
of $\yt$ on $\xt$ is given by \cc{GEWK82},
\EQN
\fyx=\intpp\fyxlam\dlamtpi \mbox{ where }
\fyxlam=ln\frac{|f_X(\lam)|}{|f_e(\lam)|}        \lab{fyxlam}
\ENN
and $f_e(\lam)$ is assembled 
(following \cc{GEWK82}) as follows.

Introduce 
$W=\sgyxe\sgxe\upm$
and note that
$\epyt-W\epxt$ is uncorrelated with $\epxt$ and has variance
$\sgepygx=\sgye-\sgyxe\sgxe^{-1}\sgxye$. Then rewrite (\ref{issd}) as
\EB
\abL{\xt}{\yt}
&=&
\aobcodL{{\hxxz}}{{\hyxz}}{{\hxyz}}{{\hyyz}}
\aobcodL{I}{W}{0}{I}
\aobcodL{I}{-W}{0}{I}
\abL{\epxt}{\epyt}\\
&=&
\aobcodL{{\hxxz+\hxyz W,}}{{\hyxz+\hyyz W,}}{{\hxyz}}{{\hyyz}}
\abL{\epxt}{\epyt-W\epxt}
\EE
This corresponds to (3.3) in \cc{GEWK82} and yields the
following expressions corresponding to those in \cc{GEWK82}.
\EBN
\fxz&=&\fez +\hxyz\sgepygx\hxy^T(\Lm) \lab{fxz}\\
\fez&=&\hex(\Lp)\sgxe \hex^T(\Lm)                             \lab{fez}\\
\hex(\Lp)&=&\hxxz+\hxyz W \nn 
\EEN
Using the SS expressions above we rewrite $\hex(L)$
in a form more suited to computation as,
\EBN
\hex(\Lp)&=&[\cx\liima B_o+I]                    \lab{hex}\\
B_o&=&B_X+B_Y\sgxye^T\sgepx\upm
=B\aob{{\sgepx}}{{\sgyxe}}\sgepx\upm \nn
\EEN
Note that then, using Theorem II,
\EBN
\fyx&=&\intpp ln|\fxlam|\dlamotpi-\intpp ln|\felam|\dlamotpi \nn\\
&=&ln|\Om_X|-ln|\sgxe|
=ln\frac{|\Om_X|}{|\sgxe|}            \lab{fyx}
\EEN
Clearly, with $L=exp(-j\lam)$, $\fxz\geq \fez\Ra \fyx\geq 0$.

Also the instantaneous causality measure is,
\EBN
\fydx&=&ln\frac{|\sgxe||\sgye|}{|\Sg_\ep|}       
=ln\frac{|\sgxe||\sgye|}{|\sgepxgy||\sgye|}
=ln\frac{|\sgxe|}{|\sgepxgy|}     \lab{fydx}\\
\sgepxgy&=&\sgxe-\sgxye\sgye^{-1}\sgyxe \nn
\EEN
Clearly $\sgxe\geq \sgepxgy$ so that $\fydx\geq 0$.

Introduce the normalised cross covariance based matrix,
$\Gam_{x,y}=\sgye^{-\half}\sgyxe\sgxe\upm\sgxye\sgye^{-\half}$.
Then using a well known partitioned matrix determinant
formula \cc{NEUD99} we find $\fydx=ln|I-\Gam_{x,y}|$.
This means that the instantaneous causality measure
depends only on the canonical correlations 
(which are the eigenvalues of $\Gam_{x,y}$) between
$\epxt,\epyt$, \cc{SEBR84},\cc{KSIH72}.

To implement these formulae,
we need  expressions for $\Om_X,\Om_Y,\fxlam$. To get them \cc{GEWK82}
fits  separate models to each of $\xt$ and $\yt$. But this causes
positivity problems with $\fyxlam$ \cc{DING06}.
Instead we obtain the required quantities from the correct submodel
obtained in the previous section. We have,\\
{\bf Theorem Va}.
The GEMs can be obtained
from the joint ISS model (\ref{issc}) 
and the submodel in Theorem II, as follows,

(a) $\fyx=ln\frac{|\Om_X|}{|\sgxe|}$ where
$\omx$ is got from the submodel in Theorem II.

(b)  The frequency domain GEM $\fyxlam$ (\ref{fyxlam}) can be computed
from (\ref{fez}),(\ref{fxlam}),(\ref{hex}).

\noi And $\Om_Y,\fxy,\fxylam$ can be obtained similarly.

Now pulling all this together with the help of result III we have
an extension of the results of \cc{GEWK82} to the
state space/VARMA case.

\noi{\bf Theorem Vb}: 
For the joint ISS model (\ref{issc}),

(a) $\fyx\geq 0,\fydx\geq 0$ and 
 $\fyx+\fydx=ln\frac{|\Om_X|}{|\sgepxgy|}$.

(b) $\yt$ dn-wgc $\xt$ iff,
$\fxz=\fez$ which holds iff $\fyx=0$ i.e. iff $\Om_X=\sgxe$.

(c) $\yt$ dn-sgc $\xt$ iff  $\fxz=\fez$ and $\sgepxgy=\sgxe$
i.e. iff  $\fyx=0$ and $\fydx=0$ i.e. 
iff $\fyx+\fydx=0$ i.e. iff $\Om_X=\sgepxgy $.

\noi{\it Remarks}. 

(i) A very nice nested hypothesis testing 
explanation of the decomposition (\ref{decomp})
is given by Parzen in the discussion to \cc{GEWK82}.

(ii) It is straightforward to see that the GEMs
are unaffected by scaling of the variables.
This is a problem for other GC measures \cc{EDIN10}.

(iii) For completeness we state extensions of the inferential
results in \cc{GEWK82} without proof.
Suppose we fit a SS model to data $\zt,t=1,\cdots,T$
using e.g. so-called state space subspace methods \cc{DEIS95},\cc{BAUR05}
or VARMA methods in e.g. \cc{LUTK93}.
Let $\fyxh,\fxyh,\fydxh,\fxoyh$ be the corresponding GEM estimators.
If we denote true values
with a superscript $0$, we find under some regularity conditions: 
\EB
&&H_0:\fyx^0=0\Ra
T\fyxh\Ra \chisq_{2n\px\py},\mbox{ as }T\ra\infty\\
&&\mbox{and } H_0:\fydx^0=0\Ra T\fydxh\Ra \chisq_{\px\py}
\EE
So to test for strong GC we put these together,
\EQ
H_0:\fyx^0=0,\fydx^0=0\Ra
T(\fyxh+\fydxh)\Ra\chisq_{(2n+1)\px\py}
\EN
Together with similar asymptotics for $\fxyh,\fxoyh $
we see that the fundamental decompositon (\ref{decomp})
has a sample version involving
a decomposition of a chi-squared into sums of 
smaller chi-squared statistics.

(iv) \cc{DING06} attempts also to derive
$\fyx$ without fitting separate models to $\xt,\yt$.
However the proposed procedure to compute $\fxlam$
involves a two sided filter and is thus in error.
The only way to get $\fxlam$ is by spectral factorization
(which produces one-sided or causal filters)
as we have done.

(v) Other kinds of causality measures have emerged
in the literature e.g. \cc{KAMI91} but it is not known whether
they obey the properties in theorems IVa,IVb. However these
properties are crucial to our subsequent analysis.

\section{\bf Effect of Filtering on Granger Causality Measures}
Now the import of the frequency
domain GEM becomes apparent since it allows us 
to determine the effect of one-sided
(or causal) filtering on GC.

We need to be clear on the situation envisaged here.
The unfiltered time series are the underlying
series of interest but we only have access to the filtered
time series. So we can only find the GEMs from
the spectrum of the filtered time series. What we need to know
is when those \dae filtered \dae GEMs are the same as the 
underlying GEMs.
We have,\\
\noi{\bf Theorem VI}.
Suppose we filter $\zt$ with a stable, full rank, one-sided filter 
$\Phi(L)=\aobcod{\phixxz}{0}{0}{\phiyyz}$ then,

(a) If $\Phi(L)$ is minimum phase then 
the GEMs (and so GC) are unaffected by filtering.

(b) If
$\Phi(L)$ has the form $\Phi(L)=\psi(L)\Phib(L)$ where $\psi(L)$
is a scalar all pass filter and $\Phib(L)$ is
stable, minimum phase
then the GEMs (and so GC) are unaffected by filtering.

(c) If $\Phi(L)$ 
is nonminimum phase and case (b) does not hold
then the GEMs (and so GC) are changed by filtering.

\noi{\it Proof}. Denote $\zbt=\Phi(L)\zt=\Phi(L)H(L)\ept$
by Result II(a). 
Then for the frequency
domain GEM we need to find, 
$\fbyxlam=ln\frac{|\fbxlam|}{|\hbex(L)\sgepx\hbex^T(\Lum)|}$
where $L=\exp(-j\lam)$.
We find trivially that. $|\fbxlam|=|\Phix(L)\fxlam\Phix(\Lum)|
=|\Phix(L)||\fxlam||\Phix(\Lum)|$. 
Finding $\hbex(L)$ is
much more complicated; we need the minimum phase vector
moving average or \stsp model corresponding to (\ref{issd}). 
Taking $\Phi(L)$ to be non-minimum
phase we carry out a spectral factorization,
$\fbzlam=\Hb(L)\Sgb\Hb^T(\Lum)$ 
where $\Hb(L)$ is causal, stable, 
 minimum phase with $\Hb(0)=I$
and then from
appendix C, $\Hb(L)$ can be written,
$\Hb(L)=\Phi(L)H(L)D(\Lum),D(\Lum)=JE^T(\Lum)\Jb\upm$
where $E(L)$ is all pass and $J,\Jb$ are 
constant matrices (Cholesky factors).
Writing this in partitioned form,
\EQ
\Hb(L)=\aobcodL{\Phix(L)}{0}{0}{\Phi_Y(L)}
\aobcodL{\hxxz}{\hyxz}{\hxyz}{\hyyz}
\aobcodL{\dxx(\Lum)}{\dyx(\Lum)}{\dxy(\Lum)}{\dyy(\Lum)}
\EN
yields $\hbex(L)=\Phix(L)\kbex(L)$ where,
\EB
\kbex(L)&=&\hxxz(\dxx(\Lum)+\dxy(\Lum)\sgepyx\sgepx\upm)\\
&+&\hxyz(\dyy(\Lum)\sgepyx\sgepx\upm+\dyx(\Lum))
\EE
Thus  in $\fbyxlam$ 
the $|\Phix(L)|$ factors cancel giving,
$\fbyxlam=ln\frac{|\fxlam|}{|\kbex(L)\sgepx\kbex^T(\Lum)|}$.
This will reduce to $\fyxlam$ iff 
$\kbex(L)=\hex(L)K(\Lum)$ 
where $K(\Lum)$ is all pass
which occurs iff
$\dxy(\Lum)=0,\dyx(\Lum)=0,
\dxx(\Lum)=\psi(\Lum)I,\dyy(\Lum)=\psi(\Lum)I$
where $\psi(\Lum)$ is a scalar all-pass filter.
Results (a),(b),(c) now follow.

We now give two examples.

{\it Example I}. Differential delay.
Suppose $\aob{\xt}{\yt}=\aobcod{1}{\rho}{0}{1}\aob{\at}{\bsubt}$
and $\Phi(L)=\aobcod{1}{0}{0}{L}$. So the two series are white noises
that exhibit an instantaneous GC. The filtering delays one series relative to the other.
Then we have,
$\zbt=\aobcod{1}{0}{0}{L}\aobcod{1}{\rho}{0}{1}\aob{\at}{\bsubt}
=\aobcod{1}{L\rho}{0}{1}\aob{\at}{\bsubt}=\Hb(L)\aob{\at}{\bsubt}$. And 
we see that $\Hb(0)=I$ while $\Hb(L)$ is stable, causal and invertible,
indeed $\Hb\upm(L)=\aobcod{1}{-L\rho}{0}{1}$. Thus we see that
the differential delay has introduced a 
spurious dynamic GC relation and
the original purely instantaneous GC is lost.

{\it Example II}. fMRI Hemodynamic 
Response is non-minimum phase.
A number of stylized or `canonical' HRFs
based on the double gamma (i.e. difference of
two gamma functions)
have been presented in the literature
e.g. \cc{HENS03},\cc{GLOV99}. 
These stylized HRFs capture two essential
features of empirical HRFs; namely
a slow rise to a peak followed by
a small negative undershoot.
And past practice
has been to use one of them for all voxels
in a slice or even volume.
Here we illustrate with a motor cortex HRF
from \cc{GLOV99} 
\EQ
h(t)=f_a(\frac{t}{\tau_am})^m\epu{-(t/\tau_a-m)}
-f_b\al(\frac{t}{\tau_bp})^p\epu{-(t/\tau_b-p)}
\EN
where $(\tau_a,m)=(1.1,5)$ and $(\tau_b,p)=(.9,12)$
while $\al=.4$. 
Also we have scaled each term to have maximum
value of $1$.
Here $f_a,f_b$ are amplitudes to be found
in a model fitting exercise.
In Fig.1 
we show a plot of the HRF with
$f_a=1=f_b$ and the zeros on a log scale.
One zero has magnitude $>1$ showing
the HRF is non-minimum phase.


\begin{figure}
\begin{center}
\resizebox{!}{8cm}{\includegraphics{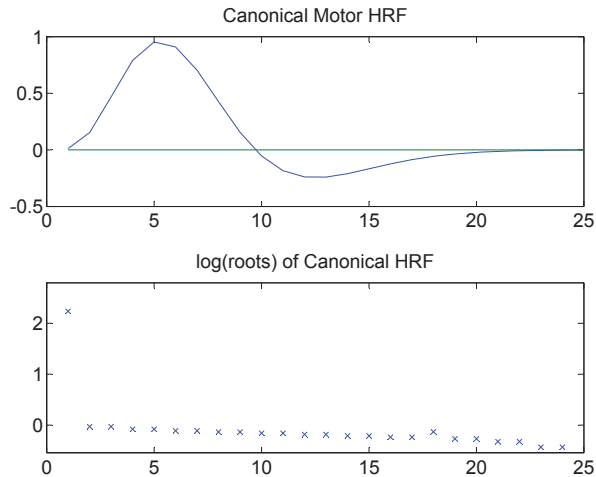}}
\caption{Canonical Motor Cortex HRF and its (log) Roots}
\end{center}
\label{hrfm}
\end{figure}

\section{\bf Downsampling and Forwards Granger Causality}
\sco
We now consider to what extent GC 
is preserved under downsampling.

Using the sampled notation of our discussion above,
and defining $\zbk=\aob{{\xbk}}{{\ybk}}$,
we have the following result:\\
\noi{\bf Theorem VII}.
Forwards Causality.

(a) If $\yt$ dn-sgc $\xt$ then 
$\ybk$ dn-sgc $\xbk$.

(b) If $\yt$ dn-wgc $\xt$ then
in general $\ybk$ wgc $\xbk$.

\noi{\it Remarks}.

(i) Part (a) is new although technically a special
case of a result of the author\dae s 
established in a non SS framework.

(ii) We might consider taking part(b)
as a formalization
of long standing folklore in econometrics
\cc{CREI87},\cc{MARC85} that downsampling
can destroy unidirectional Granger causality.
However that same folklore is flawed because
it failed to recognize the possibility of (a).
The folklore is further flawed because it failed
to recognize the more serious reverse problem
discussed below.

\noi\ul{ Proof of (a)}. 
We use the partitioned expressions
in the discussion leading up to result III.
We also refer to the discussion leading up to Theorem VI.

This allows us to write two decompositions.
Firstly $\wk=\wxk+\wyk$ where
\EQ
\wxk=\sgwm A^{m-i}\bx\ep_{X,km+i-1}\;\;,
\wyk=\sgwm A^{m-i}\by\ep_{X,km+i-1}
\EN
From result III and the definition of dn-sgc
\EQN
\wxk\mbox{ is uncorrelated with }w_{Y,l}\mbox{ for all }k,l             \lab{unca}
\ENN
The other decomposition is
$\epbk=\aob{{\epbxk}}{{\epbyk}}$
and 
\EQN
\epbxk\mbox{ is uncorrelated with }\bar{\ep}_{Y,l}\mbox{ for all }k,l      \lab{uncb}
\ENN
Next we note from Theorem IV that
$\cx(\Lm I-\aum)\upm A^p\by=\cx\sgwi A^{m+r-1}\by L^{r}=0$
for all $p\geq 0$.
Thus we deduce 
\EQN
\cx\wyk=0,\mbox{ for all }k                           \lab{wyk}
\ENN
We can now write
\EBN
\xbk&=&\cx(\Lm I-\aum)\upm\wk+\epbxk \nn\\
&=&\cx(\Lm I-\aum)\upm\wxk+\epbxk                         \lab{xbk}\\
\ybk&=&\cy\zimaumm\wxk+\cy\zimaumm\wyk+\epbyk \nn
\EEN
Based on (\ref{xbk}) we now introduce the ISS model for $\xbk$
\EQ
\xbk=\cx\zimaumm\kcx\nuxk+\nuxk
\EN
where $\nuxk$ is the innovations sequence.
Using this we introduce the estimator $\kcx\nuxk$ of $\wxk$
and the estimation error $\wtxk=\wxk-\kcx\nuxk$.
Below we show
\EQN
\wtxk \mbox{ is uncorrelated with }\nu_{X,l}\mbox{ for all }k,l   \lab{uncc}
\ENN
We thus rewrite the model for $\ybk$ as,
\EB
\ybk&=&\cy\zimaumm\kcx\nuxk+\zetak\\
\zetak&=&\cy\zimaumm[\wyk+\wtxk]+\epbyk
\EE
Now we can construct an ISS model for 
$\zetak=(I+\cy\zimaumm\kcy)\nuyk$ where $\nuyk$ is the
innovations sequence. In view of (\ref{unca},\ref{uncb},\ref{uncc})
$\nuxk$ and $\nu_{Y,l}$ are uncorrelated for all $k,l$.
Thus we have constructed the joint ISS model
\EB
\abL{{\xbk}}{{\ybk}}&=&\bar{H}(z)\abL{{\nuxk}}{{\nuyk}}\\
\bar{H}(z)&=&\aobcodL{I+\cx\zimaumm\kcx}{\cy\zimaumm\kcx}
{0}{I+\cy\zimaumm\kcy}
\EE
From this we deduce that $\ybk$ dn-sgc $\xbk$ as required.

\ul{Proof of (\ref{uncc})}. Consider then\\
$E(\wtxk\nuxl^T)=E(\wxk-\kcx\nuxk)\nuxl^T)=E(\wxk\nuxl^T-\kcx E(\nuxk\nuxl^T)$.
The second term vanishes for $k\neq l$. The first term vanishes for $k>l$
since $\wxk$ is uncorrelated with the past and hence $\nuxl$;
for $l>k$ it vanishes since $\nuxl$ is uncorrelated with the past.
For $k=l$ it vanishes by the definition of $\kcx$ \cc{KAIL00}.

\ul{Proof of (b)}.
A perusal of the proof of (a) shows that
we cannot construct the block lower triangular
joint ISS model; in general we obtain a full block ISS model.

\section{\bf Downsampling and Reverse Granger Causality}

We now come to the more serious issue of whether unidirectional
Granger causality might arise from downsampling even though not 
present on the original timescale.
To establish this we have simply to exhibit a numerical
example but that is not as simple as one might hope.

\subsection{\bf Simulation Design}
Designing a procedure to generate a wide class of
examples of spurious causality
is not as simple as one might hope.
We develop such a procedure for a 
bivariate vector autoregression of
order one; a bivariate VAR(1). On the one hand this is
about the simplest example one can consider;
on the other hand it is general enough to generate
important behaviours.

The \biva VAR(1) model is then,
\EQ
\aob{\xt}{\yt}=A\aob{x_{t-1}}{y_{t-1}}+\aob{\epxt}{\epyt}\;\;,
A=\aobcod{\phix}{\gamy}{\gamx}{\phiy}\;\;,
\Sg=\aobcod{\sg_a^2}{\rho\sg_a\sg_b}{\rho\sg_a\sg_b}{\sg_b^2}
\EN
where $\Sg$ is the \varie matrix of the zero mean white 
noise $\aob{\epxt}{\epyt}$;
$\rho$ is a correlation.

We note that this model can be written as an ISS model
with parameters, $A,I,-A,\Sg$. Hence all the computations
desribed above are easily carried out.

But the real issue is how to select the parameters.
By a straightforward scaling
argument it is easy to see the we may set $\sg_a=1=\sg_b$
without loss of generality. Thus we need to choose
only $A,\rho$.

Some reflection shows that there are two issues.
Firstly we must ensure
the process is \stty i.e. for the \evals $\lamw,\lamt$
of $A$ we must have $|\lamw|<1,|\lamt|<1$.
Secondly to design a simulation we need to choose
$\fyx,\fxy$; but these quantities depend on the parameters
$A,\rho$
in a highly nonlinear way so it is not obvious how to
do this. And five parameters is already too many
to pursue this by trial and error.

For the first issue we have 
$trace(A)=\lamw+\lamt=\phix+\phiy$
and $det(A)=\lamw\lamt=\phix\phiy-\gamx\gamy$.
Our approach is to select $\lamw,\lamt$ and then
find $\phix,\phiy$ to satisfy
$\phix+\phiy=\lamw+\lamt,\phix\phiy=\lamw\lamt+\gamx\gamy$.
This requires solution of a quadratic equation. If we
denote the solutions as $r_{+},r_{-}$ then
we get two cases: 
$(\phix,\phiy)=(r_{+},r_{-})$ and $(\phiy,\phix)=(r_{+},r_{-})$.
This leaves us to select $\gamx,\gamy$.

In Appendix B we show that 
$\fyx=ln\frac{\sg_x^2}{\sg_a^2}
\geq ln(1+\xix)$ where
$\xix=\gamx^2(1-\rho^2)$.
And similarly $\fxy\geq ln(1+\xiy)$
where $\xiy=\gamy^2(1-\rho^2)$.
But we also show that $\xix=0\Ra \fyx=0$
and $\xiy=0\Ra \fxy=0$.
So we select $\xix,\xiy$ thereby setting
a lower bounds to $\fyx,\fxy$. This seems to be the best
one can do and as we see below works quite well.
So given $\xix,\xiy$ compute
$\gamx=\pm \frac{\sqrt{\xix}}{\sqrt{1-\rho^2}}$
and $\gamy=\pm \frac{\sqrt{\xiy}}{\sqrt{1-\rho^2}}$.
This gives four cases and together with the
two cases above yields eight cases.

This is not quite the end of the story since 
the $\gamx,\gamy$ values need to be consistent with
the $\phix,\phiy$ values. Specifically the quadratic
equation to be solved for $\phix,\phiy$ must have real roots.
Thus the discriminant must be $\geq 0$. So
$(\phix+\phiy)^2-4(\phix\phiy)
=(\lamw+\lamt)^2-4(\lamw\lamt+\gamx\gamy)
=(\lamw-\lamt)^2-4\gamx\gamy\geq 0$.
There are four cases; two with real roots, two with complex roots.

If $\lamw,\lamt$ are real
then we require 
$(\lamw-\lamt)^2\geq 4\gamx\gamy
\Ra (\lamw-\lamt)^2\geq 
4sign(\gamx\gamy)\frac{\sqrt{\xix\xiy}}{1-\rho^2}$.
This always holds if $sign(\gamx\gamy)\leq 0$.
If $sign(\gamx\gamy)>0$  then we have a 
binding constraint which restricts the sizes of $\xix,\xiy$.

If $\lamw,\lamt$ are complex conjugates
then $(\lamw-\lamt)^2$ is negative.
If $sign(\gamx\gamy)\geq 0$ then the condition never holds.
If $sign(\gamx\gamy)<0$ then there is a binding
constraint which restricts the sizes of $\xix,\xiy$.
In particular note that if $sign(\gamx\gamy)=0$ then
one cannot have complex roots for $A$.
We now use this design procedure to illustrate reverse causality.

\subsection{\bf Computation}
We describe the steps used to generate the results below.
We assume the \stsp model for $\zt=\aob{\xt}{\yt}$ comes
in ISS form.
Since standard state space subspace model fitting
algorithms \cc{LARI83},\cc{OBDM96},\cc{BAUR05} generate ISS models
this is a reasonable assumption. Otherwise
we use result I to generate the corresponding
ISS model.

Given a sampling multiple $m$ we first
use Theorem III to generate the subsampled ISS model and hence
$\Sgep\upcm$. To obtain the GEMs
we use Theorem I to generate the marginal models for
$\xt,\yt$ yielding $\Om_X\upcm,\Om_Y\upcm$.
And now $\fyx\upcm,\fxy\upcm$ are gotten from
the formulae (\ref{fyx}),(\ref{fydx}) and the comment following 
Theorem Va.

\subsection{\bf Scenario Studies}
We now illustrate
the various results above with some bivariate simulations.

%
\noi\ul{Example 1}. GEMs decline gracefully.

\noi\ul{Table 1}. GEMs for various sampling intervals for\\
$(\lamw,\lamt,\xix,\xiy,\rho)
=(-.95 \expj{j\times .1}, -.95 \expj{-j\times .1},1.5,.2,.2)$
$\Ra A=\aobcod{-.204}{.452}{-1.24}{-1.69}$.\\

\bt{|l|l|l|l|l|l|l|l|l|l|l|}\hl 
m      & 1 & 2 & 3 & 4 & 5 & 6 & 10 & 20 & 30 & 40 \\\hl 
$\fyx$
&   1.3761  
&   1.657  

&   1.408   

&   1.169   

&   0.994   

&   0.864   

&   0.551  

&   0.151   

&   0.001  
&   0.014   
\\\hl
$\fxy$ 
&  0.19834

& 0.253

&   0.287

&  0.308

&  0.319

&  0.322

&  0.293

&  0.109

&  0.001

&  0.011
\\\hl 
\et\\

Here, for the underlying process,
$y$ pushes $x$ much harder than $x$ pushes $y$.
This pattern is roughly preserved with slower sampling,
but the relative strengths change.

\noi\ul{Example 2}. GEMs Reverse.

\noi\ul{Table 2}. GEMs for various sampling intervals for\\
$(\lamw,\lamt,\xix,\xiy,\rho)
=(.95 \expj{j\times .1}, .95 \expj{-j\times .1},1.5,.2,.2)$
$\Ra A=\aobcod{1.69}{.452}{-1.24}{.204}$.\\

\bt{|l|l|l|l|l|l|l|l|l|l|l|}\hl 
m      & 1 & 2 & 3 & 4 & 5 & 6 & 10 & 20 & 30 & 40 \\\hl 
$\fyx$
&   0.92983   
&   0.879   

&   0.766  
&   0.683   
&   0.62   
&   0.57  

&   0.418   

&   0.131   

&   0.001   

&   0.013  
\\\hl
$\fxy$ 
&   1.0476

&   1.824

&   2.006

&   1.795

&   1.527

&   1.3

&  0.751

&   0.18

&   0.002

&   0.016
\\\hl 
\et\\

In this case the underlying processes push eachother with
roughly equal strength. But subsampling yields a false picture
with $x$ pushing $y$ much harder than the reverse.

\noi\ul{Example 3}. Near Equal Strength Dynamics
Becomes Nearly Unidirectional.

\noi\ul{Table 3}. GEMs for various sampling intervals for\\
$(\lamw,\lamt,\xix,\xiy,\rho)
=(.995 , -.7,1,.5,.7)$
$\Ra A=\aobcod{1.45}{-.84}{1.18}{-1.16}$.\\

\bt{|l|l|l|l|l|l|l|l|l|l|l|}\hl 
m      & 1 & 2 & 3 & 4 & 5 & 6 & 10 & 20 & 30 & 40 \\\hl 
$\fyx$
&   1.487   
&   0.051  

&   0.245  
&   0.057   
&   0.111  
&   0.052 

&   0.038   

&   0.019  

&   0.012  

&   0.009  
\\\hl
$\fxy$ 
&   1.685

&   0.167

&   0.638

&   0.258

&   0.467

&   0.294

&  0.289

&   0.212

&   0.159

&   0.125
\\\hl 
\et\\

In this case the underlying relation is one of near
equal strength feedback interconnection.
But almost immediately a very unequal relation appears
under subsampling which soon decays to a near
unidirectional relation.

\noi\ul{Example 4}. Near Unidirectional Dynamics Becomes Near
Equal Strength.

\noi\ul{Table 4}. GEMs for various sampling intervals for\\
$(\lamw,\lamt,\xix,\xiy,\rho)
=(.99 \expj{j\times .25}, .95 \expj{-j\times .25},.1,3,-.8)$
$\Ra A=\aobcod{1.883}{2.236}{-0.408}{0.036}$.\\

\bt{|l|l|l|l|l|l|l|l|l|l|l|}\hl 
m      & 1 & 2 & 3 & 4 & 5 & 6 & 10 & 20 & 30 & 40 \\\hl 
$\fyx$
&   0.023   
&   0.284  

&   0.381  
&   0.428   
&   0.451 
&   0.457 

&   0.309   

&   0.454   

&   0.407   

&   0.168  
\\\hl
$\fxy$ 
&  2.937

&   2.384

&   1.617

&   1.243

&   1.019

&   0.859

&   0.372

&  0.532

&   0.453

&   0.178

\\\hl 
\et\\

In this case a near unidirectional dynamic relation
immediately becomes one of significant but unequal
strengths and then one of near equal strength.

There is nothing pathological about these examples
and using the design procedure developed above
it is easy to generate other similar kinds of examples.
They make it emphatically clear that GC cannot be 
reliably discerned from
subsampled data.
\section{\bf Conclusions}

This paper has given a theoretical and computational analysis
of the use of Granger casuality in fMRI.
There were two main issues:
the effect of downsampling and the effect of
hemodynamic convolution.
To deal with these issues a number of
novel results in multivariate time series
and Granger causality were developed 
via state space methods
as follows.

\ben
\item[(a)] Computations of submodels via the DARE (Theorems I,IV).
\item[(b)] Reliable computation of GEMs via the DARE
(Theorems Va,Vb).
\item[(c)] Effect of filtering on GEMs (Theorem VI).
In particular the destructive effect of the
non-minimum phase property of HRFs.
\item[(d)] Computation of downsampled models via the DARE.
\een

Using these results we were able to develop, in
section 8,  a framework for generating downsampling induced
spurious Granger causality 'on demand' and provided a number
of illustrations.


All this leads to the conclusion that
that Granger causality analysis of fMRI data cannot be used to discern
neuronal level driving relationships .
Not only is the time-scale too slow but
  even with faster sampling
 the non-minimum phase aspect
 of the HRF will still compromise the method.

Future work would naturally include an extension
of the Granger causality results to handle
the presence of a third vector time series.
And also extensions to deal with time-varying Granger
causality. 
Non-Gaussian versions could mitigate the non-minimum phase
problem to some extent but there does not seem
to be any evidence for the non-Gaussianity of fMRI data.
Extensions to nonlinear Granger causality
are currently of great interest but need a considerable development.

\section{\bf Stabilizability, Detectability and DARE}
In this section we restate and modify for our purposes
some standard state space results.
We rely mostly on
\cc{KAIL00}[Appendices E,C].

We denote an eigenvalue of a matrix by $\lam$
and a corresponding eigenvector by $q$.
We say $\lam$ is a \ul{stable} eigenvalue if $|\lam|<1$;
otherwise $\lam$ is an \ul{unstable} eigenvalue.

\subsection{\bf Stabilizability}
The pair $(A,B)$ is \ul{controllable} if
there exists a matrix $G$ so that $A-BG$ is stable
i.e. all eigenvalues of $A-BG$ are stable.
$(A,B)$ is controllable iff any of the following conditons hold,

(i) Controllability matrix:
$\calc=[B,AB,\cdots, A^{n-1}B]$ has rank $n$.

(ii) Rank Test:
$rank[\lam I-A, B]=n$ for all eigenvalues $\lam$ of $A$.

(iii) PBH test:
There is no left eigenvector
of $A$ that is orthogonal to $B$ i.e. if $q^TA=\lam q^T$
then $q^TB\neq 0$.

\noi The pair $(A,B)$ is \ul{stabilizable} if:
$rank[\lam I-A,B]=n$ for all \ul{unstable} eigenvalues of $A$.
Three useful tests for stabilizability are:

(i) PBH Test:
$(A,B)$ is stabilizable iff there is no left eigenvector of $A$
corresponding to an \ul{unstable} \eval
that is orthogonal to $B$ i.e. if 
$q^TA=\lam q^T$ and $|\lam|\geq 1$
then $q^TB\neq 0$.

(ii) $(A,B)$ is stabilizable if $(A,B)$ is controllable.

(iii) $(A,B)$ is stabilizable if $A$ is stable.

\subsection{\bf Theorem DARE}
\cc{KAIL00}(Theorem E6.1, Lemma 14.2.1,section 14.7)\\
Under conditions, N,St,De the DARE has a unique positive 
semi-definite solution $P$
which is stabilizing, i.e. $\As-\Ks C$ is a stable matrix.
Further if we initialize $P_0=0$ then 
$\Pt$ is nondecreasing and $\Pt\ra P$ as $t\rai$.

\noi{\it Remarks}.

(i)
$\As-\Ks C$ stable means $A-KC$ is 
stable (see below). And this implies
that $(A,K)$ is controllable (see below).

(ii) 
Since $V\geq R$ then $N\Ra V$ is positive definite.

\noi\ul{\it Proof of (i)}.
We first note (taking limits in) \cc{KAIL00}(equation 9.5.12)
$\Ks=K-SR\upm$.
We have then $\As-\Ks C=A-SR\upm C-(K-SR\upm)C=A-KC$.
So $\As-\Ks C$ is stable iff $A-KC$ is stable.
But then $(A,K)$ is controllable.

\subsection{\bf Detectability}
The pair $(A,C)$ is \ul{detectable} if $(A^T,C^T)$ is stabilizable.

\noi{\it Remarks}.

(i)
If $\As$ is stable (all eigenvalues have modulus $<1$)
then $S,D$ automatically hold.

(ii)
Condition D can be replaced with the detectability
of $(A,C)$ which is the way \cc{KAIL00} states the result.
We show equivalence below (this is also noted
in a footnote in \cc{KAIL00}(section 14.7).

\noi\ul{\it Proof of Remark(ii)}.
Suppose $(\As,C)$ is detectable but $(A,C)$ is not.
Then by the PBH test there is a right eigenvector $p$ of $A$
corresponding to an unstable \eval of $A$
with $Ap=\lam p,Cp=0$. But then $\As p=(A-SR\upm C)p=Ap=\lam p$
while $Cp=0$ which contradicts the detectbility of $(\As,C)$.
The reverse argument is much the same.\\

\ul{Proof of Result I}.
(a) follows from the discussion leading to theorem DARE.
(b) follows from the remarks after  theorem DARE.
%

\section{\bf GEMs for Bivariate VAR(1)}
Applying formula (\ref{fxz}) and reading off $\hxx $ etc. from
the VAR(1) model yields,
\EB
\fxlam&=&[(1+\phiy^2-2\phiy cos(\lam))\sga^2
+\gamx^2\sgb^2-2\rho\sga\sgb(\gamx\phiy-\gamx cos(\lam))]
/|D(\epu{j\lam})|^2\\
D(L)&=&(1-\phix L)(1-\phiy L)-\gamx\gamy L^2
\EE
This can clearly be written as an ARMA(2,1) spectrum
$\sg_x^2|1-\thtx\emu{j\lam}|^2/|D(\epu{j\lam})|^2$.
Equating coefficients gives
\EB
\gamo&=&\sg_x^2(1+\thtx^2)
=(1+\phiy^2)\sga^2+\sgb^2\gamx^2-2\rho\sga\sgb\gamx\phiy
=\sga^2(1+\xix+\dx^2)\\
\gamw&=&\sg_x^2\thtx
=\sga^2\dx
\EE
where $\xix=(1-\rho^2)\gamx^2\sgb^2$ and 
$\dx=\phiy-\rho\gamx\frac{\sgb}{\sga}$.
We thus have $\thtx^2=\gamw^2/\sg_x^4$ and using this
in the first equation gives,
$\gamo=\sg_x^2+\gamw^2/\sg_x^2$
or $\sg_x^4-\sg_x^2\gamo+\gamw^2=0$.
This has, of course, two solutions
\EB
\sg_x^2&=&\half(\gamo\pm\sqrt{\gamo^2-4\gamw^2})\\
\Ra\frac{\sg_x^2}{\sga^2}&=&
\half(\frac{\gamo}{\sga^2}\pm
\sqrt{\frac{\gamo^2}{\sga^4}-4\frac{\gamw^2}{\sga^4}})
=\half(1+\xix+\dx^2\pm\sqrt{(1+\xix+\dx^2)^2-4\dx^2})
\EE
Note that if $\xix=0$ this delivers
$\frac{\sg_x^2}{\sga^2}=\half(1+\dx^2+\sqrt{(1-\dx^2)^2}=1$.

We must choose the solution which ensures 
$|\thtx|<1
\equiv \frac{\gamw^2}{\sg_x^4}<1
\equiv \frac{\gamo}{\sg_x^4}<2
\equiv \gamo/\sga^2<2\sg_x^2/\sga^2$
$\equiv 1+\xix+\dx^2<1+\xix+\dx^2\pm \sqrt{(1+\xix+\dx^2)^2-4\dx^2}$.
And so we must choose the \dae$+$\dae solution.
Continuing, we now claim,
\EQ
\frac{\sg_x^2}{\sga^2}\geq\half(1+\xix+\dx^2+\sqrt{(1+\xix-\dx^2)^2}
=\half(1+\xix+\dx^2+1+\xix-\dx^2)
=1+\xix
\EN
This follows if,
$(1+\xix+\dx^2)^2-4\dx^2\geq(1+\xix-\dx^2)^2
\equiv (1+\xix+\dx^2)^2-(1+\xix-\dx^2)^2\geq4\dx^2
\equiv 2(1+\xix)2\dx^2\geq4\dx^2
\equiv 1+\xix\geq1$,
which holds.

\section{\bf Spectral Factorization}
Suppose $\ept$ is a white noise sequence with 
$E(\ept)=0,var(\ept)=\Sg$.
Let $G(L)$
be a stable causal possibly non-minimum phase filter.
Then $\zbt=G(L)\ept$ has spectrum
$\fbzlam=G(L)V G^T(\Lum)$ where $L=exp(-j\lam)$.
We can then find a unique causal, stable minimum phase spectral
factorization, $\fbzlam=\Go(L)V_o\Go^T(\Lum)$.
Let $V,V_o$ have Cholesky factorizations, 
$V=JJ^T,V_o=\Jo\Jo^T$
and set $\Gc(L)=G(L)J,\Goc(L)=\Go(L)\Jo$.
Then $\fbzlam=\Gc(L)\Gc^T(\Lum)=\Goc(L)\Goc^T(\Lum)$.
Since $\Goc(L)$ is minimum phase we can introduce
the causal filter $E(L)=\Goc\upm(L)\Gc(L)\Ra
E(L)E^T(\Lum)=I$.
Such a filter is called an \ul{all} \ul{pass} filter
\cc{HDE88},\cc{GREE88}. 
Now, $\Gc(L)=\Goc(L) E(L)$
or $G(L)=\Go(L)\Jo E(L)J\upm$ i.e. a
decomposition of a non-minimum phase
(matrix) filter into a product of a minimum phase
filter and an all pass filter. We can also write this as,
$\Go(L)=G(L)J E\upm(L)\Jo\upm
=G(L)JE^T(\Lum)\Jo\upm$ showing how
the non-minimum phase filter is transformed to yield a spectral factor.

%
          
\end{document}